\documentclass{gtart}

\def\ifplaintex{\expandafter\ifx\csname documentclass\endcsname\relax}

\expandafter\ifx\csname epsfbox\endcsname\relax\input epsf\fi

\ifplaintex 
\else
\fi

\def\gt{{\mathsurround=0pt\it $\cal G\mskip-2mu$eometry \&\ 
$\cal T\!\!$opology}}        

\def\gtp{{\mathsurround=0pt\it $\cal G\mskip-2mu$eometry \&\ 
$\cal T\!\!$opology $\cal P\!$ublications}}  

\def\lognumber#1{\def\thelognumber{#1}}
\def\volumenumber#1{\def\thevolumenumber{#1}}
\def\papernumber#1{\def\thepapernumber{#1}}
\def\volumeyear#1{\def\thevolumeyear{#1}}

\def\pagenumbers#1#2{\def\startpage{#1}\def\finishpage{#2}}
\def\published#1{\def\publishdate{#1}}
\def\proposed#1{\def\theproposer{#1}}
\def\seconded#1{\def\theseconders{#1}}
\def\received#1{\def\receiveddate{#1}}
\def\revised#1{\def\reviseddate{#1}}
\def\accepted#1{\def\accepteddate{#1}}

\long\def\asciiabstract#1{\long\def\theasciiabstract{#1}}

\let\\\par\let\thevolumenumber\relax\let\thepapernumber\relax
\let\thevolumeyear\relax\let\thesamplenumber\relax\let\startpage\relax
\let\finishpage\relax\let\publishdate\relax\let\receiveddate\relax
\let\reviseddate\relax\let\accepteddate\relax\let\theasciititle\relax
\let\theasciiauthors\relax
\let\theasciiabstract\relax
\let\theasciiemail\relax\let\theshortauthors\relax\let\theshorttitle\relax

\long\def\maketitlep{   

\count0=\startpage

\gt\hfill      
\hbox to 77pt{\vbox to 0pt{\vglue -15pt\epsfbox{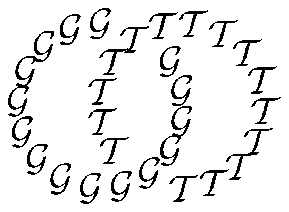}\vss}\hss}
\break
{\small\ifx\thesamplenumber\relax 
Volume \else Sample
\fi\thevolumenumber\ (\thevolumeyear)
\startpage--\finishpage\nl
Published: \publishdate}
\vglue 0.5truein plus 0.4fil minus 0.1truein

{\parskip=0pt\leftskip 0pt plus 1fil\def\\{\par\smallskip}{\ifplaintex\large
\else\Large\fi\bf\thetitle}\par\medskip}   

\vglue 0pt plus 0.1fil 

{\parskip=0pt\leftskip 0pt plus 1fil\def\\{\par}{\sc\theauthors}
\par\medskip}

\vglue 0pt plus 0.1fil 

{\small\parskip=0pt\let\newline\\
{\leftskip 0pt plus 1fil\def\\{\par}{\sl\theaddress}\par}
\expandafter\ifx\theemail\relax    
\relax\else\vglue 5pt plus 0.02fil minus 2pt\def\\{\stdspace{\rm 
and}\stdspace} 
\cl{Email:\stdspace\tt\theemail}\fi
\ifx\theurl\relax                  
\relax\else\vglue 5pt plus 0.02fil minus 2pt\def\\{\stdspace{\rm 
and}\stdspace}
\cl{URL:\stdspace\tt\theurl}\fi\par}

\vglue 7pt plus 0.3fil minus 3pt

{\bf Abstract}
\vglue 5pt plus 0.1fil minus 2pt

\theabstract

\vglue 7pt plus 0.3fil minus 3pt

{\bf AMS Classification numbers}\quad Primary:\quad \theprimaryclass

Secondary:\quad \thesecondaryclass

\vglue 5pt plus 0.3fil minus 2pt

{\bf Keywords:}\quad \thekeywords

\vglue 10pt plus 0.5fil minus 5pt

{\small  Proposed: \theproposer\hfill Received: \receiveddate\nl
Seconded: \theseconders\hfill 
\ifx\reviseddate\relax                         
Accepted: \accepteddate                        
\else
Revised: \reviseddate                          
\fi}
\eject
}       

\let\maketitlepage\maketitlep
\let\maketitle\maketitlepage

\font\phead=cmsl9 scaled 950
\font=cmsl9 scaled 1050
\font\pnum=cmbx10 scaled 913
\font\lnum=cmbx10 
\font\pfoot=cmsl9 scaled 950
\font=cmsl9 scaled 1050
\ifplaintex
\headline{\vbox to 0pt{\vskip -4.5mm\line{\small\phead\ifnum
\count0=\startpage ISSN 1364-0380 (on line)
1465-3060 (printed) \hfill {\pnum\folio}\else\ifodd\count0\def\\{ }%
\ifx\theshorttitle\relax\thetitle\else\theshorttitle\fi\hfill{\pnum\folio}
\else\def\\{ and }{\pnum\folio}\hfill\ifx\theshortauthors\relax\theauthors
\else\theshortauthors\fi\fi\fi}\vss}}
\footline{\vbox to 0pt{\vglue 0mm\line{\small\pfoot\ifnum\count0=\startpage
\copyright\ \gtp\hfill\else
\gt, Volume \thevolumenumber\ (\thevolumeyear)\hfill\fi}\vss
}}
\else
\makeatletter
\def\@oddhead{{\small\lhead\ifnum\count0=\startpage ISSN 1364-0380 (on line)
1465-3060 (printed) \hfill {\lnum\number\count0}\else\ifodd\count0
\def\\{ }\ifx\theshorttitle\relax \thetitle \else\theshorttitle\fi\hfill
{\lnum\number\count0}\else\def\\{ and }{\lnum\number\count0}
\hfill\ifx\theshortauthors\relax 
\theauthors\else\theshortauthors\fi\fi\fi}}\def\@evenhead{@oddhead}
\def\@oddfoot{\small\lfoot\ifnum\count0=\startpage\copyright\ \gtp\hfill\else
\gt, Volume \thevolumenumber\ (\thevolumeyear)\hfill\fi}
\def\@evenfoot{@oddfoot}
\makeatother
\fi

\newwrite\gtoutfile
\long\gdef\makeheadfile{  
{\def\\{, }\def\s{ }
\immediate\openout\gtoutfile head.xxx
\immediate\write\gtoutfile{To: math@arxiv.org}
\immediate\write\gtoutfile{Subject: put OR rep NNNNN:pppp}
\immediate\write\gtoutfile{--text follows this line--}
\immediate\write\gtoutfile{Proxy-for: \ifx\theasciiauthors\relax
\theauthors\else\theasciiauthors\fi\s<\ifx\theasciiemail\relax\theemail\else\theasciiemail\fi>}
\immediate\write\gtoutfile{\noexpand\\}
\immediate\write\gtoutfile{Authors: \ifx\theasciiauthors\relax
\theauthors\else\theasciiauthors\fi}
{\def\\{ }\immediate\write\gtoutfile{Title: \ifx\theasciititle\relax
\thetitle\else\theasciititle\fi}}
\immediate\write\gtoutfile{Subj-class: GT or GR or SG or ...}
\immediate\write\gtoutfile{MSC-class: \theprimaryclass\ifx\thesecondaryclass\relax\else, \thesecondaryclass\fi}
\immediate\write\gtoutfile{Journal-ref: Geom. Topol. \thevolumenumber\s
(\thevolumeyear) \startpage-\finishpage}
\immediate\write\gtoutfile{Comments: Published in Geometry and Topology at}
\immediate\write\gtoutfile{    http://www.maths.warwick.ac.uk/gt/GTVol\thevolumenumber/paper\thepapernumber.abs.html}
\immediate\write\gtoutfile{\noexpand\\}
\immediate\write\gtoutfile{}
\ifx\theasciiabstract\relax
\immediate\write\gtoutfile{\theabstract}\else
\immediate\write\gtoutfile{\theasciiabstract}\fi
\immediate\write\gtoutfile{}
\immediate\write\gtoutfile{\noexpand\\}
\immediate\write\gtoutfile{}
\immediate\closeout\gtoutfile}}  

\def\maketitlepage{\maketitlep\makeheadfile}
\let\maketitle\maketitlepage

\lognumber{232}
\volumenumber{7}\papernumber{3}\volumeyear{2003}
\pagenumbers{91}{153}
\received{1 February 2002}
\revised{18 October 2003}
\accepted{19 November 2002}
\published{17 February 2003}
\proposed{Haynes Miller}
\seconded{Ralph Cohen, Gunnar Carlsson}

\usepackage{amsmath,amssymb,amscd}

\usepackage[curve,cmtip,matrix,arrow]{xy}

\def\S{section }

\numberwithin{equation}{section}
\newtheorem{Theorem}[equation]{Theorem}
\newtheorem{Conjecture}[equation]{Conjecture}
\newtheorem{Lemma}[equation]{Lemma}
\newtheorem{Proposition}[equation]{Proposition}
\newtheorem{Corollary}[equation]{Corollary}

\theoremstyle{definition}

\newtheorem{Definition}[equation]{Definition}

\theoremstyle{remark}

\newtheorem{Remark}[equation]{Remark}
\newtheorem{Example}[equation]{Example}

\newcommand{\eqdef}{\overset{\text{def}}{=}}
\newcommand{\psb}[1]{[ \! [#1] \! ]}
\newcommand{\lsb}[1]{( \! (#1) \! )}
\newcommand{\tensor}[1]{\underset{#1}{\otimes}}

\newcommand{\slot}{\,-\,}

\newcommand{\CP}{\mathbb{C} \! P^{\infty}}

\newcommand{\Z}{\mathbb{Z}}
\newcommand{\Q}{\mathbb{Q}}

\newcommand{\C}{\mathbb{C}}

\newcommand{\xra}{\xrightarrow}

\newcommand{\Ad}{c}

\newcommand{\aff}{\mathbb{A}^{1}}
\newcommand{\affan}{\mathbb{A}^{1}_{\mathrm{an}}}
\newcommand{\affanfml}{(\mathbb{A}^{1}_{\mathrm{an}})^{\wedge}_{0}}
\newcommand{\affanstalk}{\mathbb{A}^{1}_{\mathrm{an},0}}
\newcommand{\ac}[1]{A (#1)}
\newcommand{\anomaly}[1]{\mathcal{A} (#1)}
\DeclareMathOperator{\Aut}{Aut}

\newcommand{\bara}{\overline{a}}
\newcommand{\barm}{\overline{m}}
\newcommand{\barU}{\overline{U}}
\newcommand{\boeight}{BO\langle 8 \rangle}
\newcommand{\borel}[1]{#1_{\T}}

\newcommand{\cat}[1]{\mathcal{#1}}
\newcommand{\CategoryOf}[1]{(\text{#1})}

\newcommand{\chars}{\hat{T}}
\renewcommand{\circle}{\mathbb{T}}
\newcommand{\cochars}{\check{T}}
\DeclareMathOperator{\colim}{colim}
\newcommand{\compact}[1]{#1^{+}}

\newcommand{\CQ}{C}
\newcommand{\CQhat}{\fmlgpof{C}}

\newcommand{\dfcc}{\xi}

\newcommand{\Egkv}[2]{GKV_{#1} (#2)}
\newcommand{\Esym}{E}
\newcommand{\E}[1]{\spaceof{\Esym} (#1)}
\newcommand{\EDsym}{\mathcal{E}}
\newcommand{\ED}[1]{\EDsym (#1)}
\newcommand{\EPsym}{\hat{\Esym}}
\newcommand{\EP}[1]{\spaceof{\EPsym} (#1)}
\newcommand{\ETsym}{E_{\T}}
\newcommand{\ET}[1]{\spaceof{\ETsym} (#1)}
\newcommand{\ETS}[1]{\ETsym (#1)}
\newcommand{\ETa}[2]{\spaceof{E_{\T,#2}} (#1)}
\newcommand{\ETSa}[2]{E_{\T,#2} (#1)}
\newcommand{\etorsa}{\alpha}
\newcommand{\etorsap}{\beta}

\newcommand{\fmlgpof}[1]{\widehat{#1}}

\newcommand{\Ga}{\mathbb{G}_{a}}
\newcommand{\Gah}{\widehat{\mathbb{G}}_{a}}
\newcommand{\Gaan}{\mathbb{G}_{a}^{\mathrm{an}}}
\newcommand{\gkv}{\text{ideal}}
\newcommand{\Gman}{\mathbb{G}_{m}^{\mathrm{an}}}
\newcommand{\Gm}{\mathbb{G}_{m}}
\newcommand{\GpOf}[1]{P_{#1}}
\newcommand{\glue}{\psi}

\newcommand{\halfplane}{\mathfrak{H}}
\newcommand{\h}{\halfplane}
\newcommand{\HCsym}{H}
\newcommand{\HC}[1]{\spaceof{\HCsym} (#1)}
\newcommand{\HHsym}{\mathcal{H}}
\newcommand{\HHstalk}[1]{\HHsym (#1)_{0}}
\newcommand{\HPsym}{\hat{H}}

\newcommand{\HPZsym}{HP\Z}

\newcommand{\HPCsym}{HP}

\newcommand{\HQsym}{H\Q}
\newcommand{\HQ}[1]{\spaceof{\HQsym} (#1)}
\newcommand{\HZsym}{H\Z}
\newcommand{\HZ}[1]{\spaceof{\HZsym} (#1)}

\newcommand{\HHol}[2]{\mathcal{H} (#1;#2)}

\newcommand{\iso}{\cong}
\newcommand{\I}{\mathcal{I}}
\newcommand{\Ihatsym}{\hat{I}}
\newcommand{\Ihat}[1]{\Ihatsym (#1)}

\newcommand{\Loo}{\mathcal{L}}
\newcommand{\Loops}[1]{L #1}

\newcommand{\musix}{MU\langle 6 \rangle}
\newcommand{\moeight}{MO\langle 8 \rangle}
\renewcommand{\O}[1]{\mathcal{O}_{#1}}

\newcommand{\norm}[1]{|#1|}
\newcommand{\oneforms}{\uln{\omega}}

\newcommand{\orient}{\epsilon}
\newcommand{\OverS}{\mathcal{S}}
\newcommand{\PB}{Q}
\newcommand{\point}{\ast}
\newcommand{\pr}{\wp}

\newcommand{\qfmsym}{\phi}
\newcommand{\qfm}[1]{\qfmsym (#1)}

\newcommand{\rank}{r}
\newcommand{\restr}[1]{|_{#1}}
\newcommand{\ringedspaces}{\mathcal{R}}

\newcommand{\sg}{\sigma}
\DeclareMathOperator{\spec}{spec}
\DeclareMathOperator{\spf}{spf}

\DeclareMathOperator{\spaceop}{sp}
\newcommand{\spaceof}[1]{\spaceop #1}

\newcommand{\spin}{Spin}
\newcommand{\stalk}[1]{{#1}_{0}}
\newcommand{\suchthat}{\;|\;}
\newcommand{\SubT}{Sub (\T)}

\newcommand{\uhalf}{u^{\frac{1}{2}}}
\newcommand{\umhalf}{u^{-\frac{1}{2}}}

\newcommand{\trans}{\tau}
\newcommand{\T}{\circle}

\newcommand{\spfcohborel}[2]{\spaceof{#2}(\borel{#1})}
\newcommand{\TE}[2]{\spfcohborel{#1}{\Esym}\restr{#2}}
\newcommand{\TEP}[1]{\spfcohborel{#1}{\EPsym}}
\renewcommand{\TH}[1]{\spaceof{\mathcal{H} (#1)}}
\newcommand{\THstalk}[1]{\spfcohborel{#1}{H}_{0}}
\newcommand{\THC}[1]{\spfcohborel{#1}{\HCsym}}
\newcommand{\THZ}[1]{\spfcohborel{#1}{\HZsym}}
\newcommand{\THPC}[1]{\spfcohborel{#1}{\HPCsym}}

\newcommand{\uln}[1]{\underline{#1}}

\newcommand{\VG}[1]{\cochars \otimes #1}
\newcommand{\VGW}[1]{(\VG{#1})/W}

\begin{document}
\title[Analytic equivariant sigma orientation]{The sigma orientation
for analytic\\circle-equivariant elliptic cohomology} 

\author{Matthew Ando}
\email{mando@math.uiuc.edu}

\address{Department of Mathematics\\University of Illinois at 
Urbana-Champaign\\Urbana IL 61801, USA}

\begin{abstract}
We construct a \emph{canonical} Thom isomorphism in Grojnowski's
equivariant elliptic cohomology, for virtual $\T$-oriented
$\T$-equivariant spin bundles  with vanishing Borel-equivariant second
Chern class, which is natural under pull-back of vector bundles and
exponential under Whitney sum. It extends in the complex-analytic case the
non-equivariant sigma orientation of Hopkins, Strickland, and the
author. The construction relates the sigma orientation to the
representation theory of loop groups and Looijenga's weighted
projective space, and sheds light even on the non-equivariant
case. Rigidity theorems of Witten-Bott-Taubes including
generalizations by Kefeng Liu follow.  
\end{abstract}

\asciiabstract{We construct a canonical Thom isomorphism in
Grojnowski's equivariant elliptic cohomology, for virtual T-oriented
T-equivariant spin bundles with vanishing Borel-equivariant second
Chern class, which is natural under pull-back of vector bundles and
exponential under Whitney sum. It extends in the complex-analytic case
the non-equivariant sigma orientation of Hopkins, Strickland, and the
author. The construction relates the sigma orientation to the
representation theory of loop groups and Looijenga's weighted
projective space, and sheds light even on the non-equivariant
case. Rigidity theorems of Witten-Bott-Taubes including
generalizations by Kefeng Liu follow.}

\keywords{Sigma orientation, equivariant elliptic cohomolgy, rigidity}

\primaryclass{55N34}

\secondaryclass{55N22, 57R91}
\maketitlepage


\section{Introduction}\label{sec-intro}


Let $E$ be an even periodic, homotopy commutative ring spectrum, let $C$
be an elliptic curve over $S_{E} = \spec \pi_{0}E$, and let $t$ be an
isomorphism of formal groups 
\[
   t:  \fmlgpof{C} \cong \spf E^{0} (\CP),
\]
so that $(E,C,t)$ is an elliptic spectrum in the sense of
\cite{ho:icm,AHS:ESWGTC}.   In \cite{AHS:ESWGTC}, Hopkins, Strickland,
and the author construct a canonical map of homotopy commutative ring
spectra 
\[
   \sigma (E,C,t): \musix \xra{} E
\]
called the \emph{sigma orientation}; it is conjectured in
\cite{ho:icm} that this map is the restriction to $\musix$ of a
similar map $\moeight\to E.$

Let $\T$ be the
circle group.  We expect that there is an equivariant elliptic
spectrum $(\ETsym,C,t)$ extending the nonequivariant elliptic spectrum
$(E,C,t)$, and that the 
sigma orientation extends to a multiplicative map of $\T$-equivariant spectra
\begin{equation} \label{eq:11}
\sigma_{\T} (\ETsym,C,t): (\text{$\T$-equivariant $\moeight$}) 
\xrightarrow{}
\ETsym.
\end{equation}
Note however
that the construction of $\sigma_{\T}$ requires us among other
things to say what $\T$-spectra  we have in mind for the domain and codomain.  

If $V$ is a virtual vector bundle over a space $X$, let $X^{V}$ denote
its Thom spectrum.  In principle, giving a map $\sigma_{\T}$ as in
\eqref{eq:11} should be  equivalent to specifying, for each virtual
$\T$-$\boeight$ vector bundle  (whatever that means) a trivialization
$\gamma (V)$ of $\ETS{X^{V}}$ as an $\ETS{X}$-module; the
trivialization should be stable, exponential, and natural as $V/X$
ranges over the virtual $\T$-$\boeight$ vector bundles.

The trivializations $\gamma (V)$ should be compatible with 
the nonequivariant sigma orientation in the following way.  If $X$ is
a $\T$-space, then $\borel{X}$ will denote the  Borel construction 
$E\T \times_{\T} X;$ if $V$ is a (virtual) $\T$-vector bundle over
$X$, then $\borel{V}$ will denote the corresponding (virtual) bundle
over $\borel{X}$.  A $\T$-$\boeight$ structure on a $\T$-bundle
should at least give a $\boeight$ structure to $\borel{V}$.  One
expects that the equivariant extension $\ETsym$ of an elliptic
cohomology theory $E$ comes with a 
completion isomorphism 
\begin{equation} \label{eq:39}
     \ETS{X}^{\wedge} \cong E (\borel{X}),
\end{equation}
and in particular
\begin{equation} \label{eq:35}
     \ETS{X^{V}}^{\wedge}\cong E ((\borel{X})^{\borel{V}}).
\end{equation}
The desired compatibility is that this isomorphism carries $\gamma
(V)$ to the sigma orientation  of $\borel{V}$. 

In this paper we take $\ETsym$ to be the 
complex-analytic equivariant elliptic cohomology of Grojnowski.
Let
$\Lambda\subset \C $ be a lattice in the complex plane, and let $C$
be the analytic variety $C=\C /\Lambda$.  Grojnowski
constructs a functor $\ETsym$ from 
finite $\T$-CW complexes to sheaves of $\Z/2$-graded
$\O{C}$-algebras, equipped with a natural isomorphism~\eqref{eq:39},
where the left side denotes the completion of the stalk of $\ETS{X}$ at the origin of $C$ 
(\cite{Grojnowski:Ell}; for a published account see 
\cite{Rosu:Rigidity,AndoBasterra:WGEEC}).  

The bundles for which we construct trivializations $\gamma (V)$ are
the 
\emph{virtual $\T$-oriented equivariant spin bundles $V$ with $c_{2}
(\borel{V})=0$}.  This requires some explanation.

Let $V$ be a $\T$-vector bundle over $X$.  A \emph{$\T$-orientation}
$\orient$ on $V$ is a choice of orientation $\orient (V^{A})$ on the
fixed sub-bundle 
$V^{A}$ for each closed subgroup $A$ of $\T$.   We say that $V$ is
\emph{$\T$-orientable} if it admits a $\T$-orientation; a
\emph{$\T$-oriented vector bundle} is a $\T$-vector bundle equipped
with a $\T$-orientation.  An isomorphism of
$\T$-oriented vector bundles is an isomorphism of $\T$-vector bundles
which preserves the orientations on each of the fixed sub-bundles. 

A \emph{$\T$-equivariant spin bundle} is a spin vector bundle $V$,
equipped with an action of $\T$ on its principal spin bundle $Q \to X$.  If $V$ is a $\T$-equivariant spin bundle with principal spin bundle
$Q$, then the Borel construction $\borel{Q}$ is a principal
spin bundle over $\borel{X}$; the associated spin vector bundle is
$\borel{V}$. 

Since we produce a section of the sheaf which corresponds in
Grojnowski's $\ETsym$  to the cohomology of the Thom
spectrum of a virtual bundle, it is important to identify this
sheaf precisely, rather than merely its isomorphism class.  Our
approach is
to construct the section for a pair $V= V/X = (V_{0},V_{1})$ of
$\T$-oriented equivariant spin vector bundles over $X$, representing
the virtual bundle $V_{0} - V_{1}$, and then state explicitly how  the
construction depends on the pair.  We choose not 
to write $V=V_{0} - V_{1}$, to 
emphasize that we are working with genuine vector bundles and not
mere isomorphism classes. In particular the virtual bundles over $X$
do not form an abelian group.  On the other hand, if $V$ is such a
pair, then $X^{V}$ will denote the Thom spectrum of the virtual bundle
$V_{0} - V_{1}$.

If $V= (V_{0},V_{1})$ is a pair of vector bundles over $X$, and $f: Y\to X$ is
a map, then we write $f^{*}V$ for the pair $(f^{*}V_{0},f^{*}V_{1})$
of vector bundles over $Y$.  A map
\[
   g/f: W/Y \to  V/X
\]
of pairs of vector bundles will mean a map of spaces
\[
    f: Y\to X
\]
and maps of vector bundles 
\begin{equation} \label{eq:48}
\begin{CD}
W_{i} @> g_{i} >> V_{i} \\
@VVV @VVV \\
Y @> f >> X
\end{CD}      
\end{equation}
covering $f$ for $i=0,1$.   Such a map will be called a
\emph{pull-back} if the diagram \eqref{eq:48} is a pull-back for each
$i$.  If $V$ is a pair of $\T$-vector bundles, then we write
$\borel{V}$ for the pair $(\borel{( V_{0})}, \borel{( V_{1})}).$
If $V'= (V_{0}', V_{1}')$ is another
pair of vector bundles over $X$, then we write $V\oplus V'$ for the
pair $(V_{0}\oplus V_{0}', V_{1}\oplus V_{1}')$ over $X$.

A spin vector bundle $V$ has an integral characteristic class, twice
which is the first Pontrjagin class.  It is called $\frac{p_{1}}{2}$ in
\cite{BottTaubes:Rig} and $\lambda$ in \cite{fw:asd}.  The
isomorphism  
\[
   H^{4} (BSU;\Z) \cong H^{4} (BSpin;\Z)
\]
identifies it with the second Chern class, so we call it $c_{2}$
as in \cite{AHS:ESWGTC}.
If $V= (V_{0}, V_{1})$ is a pair of spin
bundles, then we define $c_{2} (V) = c_{2} (V_{0}) - c_{2} (V_{1})$.  

If $V$ is a $\T$-equivariant vector bundle over $X$, then we write
$\ETS{V}$ for the reduced equivariant elliptic cohomology 
of its Thom space.  It is an $\ETS{X}$-module, and
if $V$ is $\T$-orientable, then it is an
invertible $\ETS{X}$-module (see
\S\ref{sec:ellipt-cohom-thom-spaces} and
\cite{Rosu:Rigidity,AndoBasterra:WGEEC}).  If $V$ and $V'$ are two
such bundles then there is a canonical isomorphism of $\ETS{X}$-modules 
\begin{equation} \label{eq:16}
    \ETS{V \oplus V'} \cong 
    \ETS{V}\tensor{\ETS{X}} \ETS{V'}.
\end{equation}
If $V$ is a $\T$-oriented vector bundle, and $f: Y\to X$ is an
equivariant map, then there is a canonical isomorphism 
\begin{equation} \label{eq:46}
   \ETS{f^{*}V} \cong f^{*} \ETS{V}
\end{equation}
of sheaves of $\ETS{Y}$-modules.

If $V= (V_{0},V_{1})$ is a pair of $\T$-oriented
vector bundles over $X$, then we define 
\[
     \ETS{V} \eqdef \ETS{V_{0}}\tensor{\ETS{X}} \ETS{V_{1}}^{-1}.
\]
Equations~\eqref{eq:16} and~\eqref{eq:46} imply 
that the assignment $V\mapsto \ETS{V}$ is
\begin{enumerate}
\item \emph{stable} in the sense that if $W$ is a $\T$-oriented vector
bundle over $X$, and if $V'= (V_{0}\oplus W, V_{1}\oplus W)$, then
there is a canonical isomorphism   
\begin{equation} \label{eq:44}
      \ETS{V'} \cong \ETS{V};
\end{equation}
\item \emph{natural} in the sense that if 
\[
g/f: W/Y \rightarrow V/X
\]
is a pull-back of pairs of $\T$-oriented vector bundles,
then there is a canonical isomorphism 
\begin{equation} \label{eq:49}
    \ETS{g/f}: \ETS{f}^{*} \ETS{V} \cong \ETS{W}
\end{equation}
of $\ETS{Y}$-algebras.
\item \emph{exponential} in the sense that if $V'= (V_{0}',V_{1}')$ is
another $\T$-oriented vector bundle, then there is a canonical isomorphism 
\begin{equation}  \label{eq:45}
    \ETS{V \oplus V'} \cong 
    \ETS{V}\tensor{\ETS{X}} \ETS{V'}
\end{equation}
extending \eqref{eq:16}.
\end{enumerate}
Our main result is the following.

\begin{Theorem}[\ref{t-th-sigma-orientation}] \label{t-th-sigma-intro}
Let $V= (V_{0},V_{1})$ be a  pair of $\T$-oriented equivariant spin
vector bundles over a finite $\T$-CW complex $X$, with the
property that $c_{2} (\borel{V})  = 0.$
Then there is a 
canonical trivialization $\gamma (V)$ of $\ETS{V}$, whose 
value in $\ETS{V}^{\wedge} \cong E (\borel{X}^{\borel{V}})$ is
the Thom class provided by the (nonequivariant) sigma orientation.

The association $V\mapsto \gamma (V)$ is stable in the sense that if
$W$ is a $\T$-oriented $\T$-equivariant spin bundle, and 
$V' = (V_{0}\oplus W, V_{1}\oplus W)$, then 
\[
    \gamma (V') = \gamma (V)
\]
under the isomorphism~\eqref{eq:44}. 
It is natural in the sense that if 
\[
  g/f: W/Y \rightarrow V/X
\]
is a pull-back of pairs of $\T$-oriented equivariant spin
bundles with vanishing Borel second Chern class, then
\[
    \ETS{f}^{*} \gamma (V) = \gamma (V')
\]
under the isomorphism~\eqref{eq:49}.
It is exponential 
in the sense that if $V'= (V_{0}',V_{1}')$ is another pair of $\T$-oriented
$\T$-equivariant spin bundles with $c_{2} (\borel{V'}) = 0$, then 
\[
   \gamma (V \oplus  V') = \gamma (V) \otimes \gamma (V')
\]
under the isomorphism~\eqref{eq:45}.  
\end{Theorem}

\begin{Remark} 
It is a result of Bott and Samelson (\cite{BS:atmss}; see
\cite{BottTaubes:Rig} or Lemma \ref{t-le-Z-connected-oriented})
that a $\T$-equivariant spin bundle is 
$\T$-orientable.  
\end{Remark}

\begin{Remark}
The class $\gamma (V)$ depends on the 
$\T$-orientation and spin structure, but otherwise requires only that
$c_{2} (\borel{V})=0$, i.e. that $\borel{V}$ \emph{admits} a
$\boeight$-structure.  We suspect that the failure of our construction
to depend on the choice of $\boeight$ structure on $\borel{V}$
reflects the fact that 
our elliptic curve $C$ comes as a quotient
$\C/\Lambda$.    For example in Lemma \ref{t-F-bara-weil}, it is
important to be able to name the point $q^{\frac{1}{n}}$.  This is
similar to the situation of the Ochanine genus, whose rigidity 
requires only that $V$ is a spin bundle.
\end{Remark}

In \cite{AndoBasterra:WGEEC}, Maria Basterra and the author
showed that, under the hypotheses of the Theorem, 
there is a global  section $\gamma (V)$ of $\ETS{V}$, whose value in
the stalk at the origin of $C$ is the sigma
orientation of $\borel{V}$; earlier Rosu \cite{Rosu:Rigidity} did the same for
the orientation associated to the Euler formal group law.  However,
those papers do not address the trivialization, naturality, and
exponential properties of the classes they construct.  

The naturality is particularly hard to discern.  
Both of the papers
\cite{Rosu:Rigidity,AndoBasterra:WGEEC} closely follow 
\cite{BottTaubes:Rig} in their implementation of the ``transfer''
argument, and so all three papers depend on meticulous choices for
integer representatives of the characters of the action of $\T$ on
$V\restr{X^{\T}}$ and of $\T[n]$ on $V\restr{\T[n]}$ and for the
orientations of $V^{\T}$ and $V^{\T[n]}$, along with surprising and
not particularly intuitive results concerning the compatibility of
these choices.   In this paper we 
show that \emph{any} choice will do, with \emph{no} effect on the
resulting Thom class, which is completely determined by the
equivariant spin bundle $V$, together with the orientations
of the bundles $V^{\T[n]}$ and $V^{\T}$. The argument is indifferent
to the parity of $n$. 

From a practical point of view there are two important new ingredients.  The
first is a careful account of the choice of representatives of
characters (rotation numbers) for the action of $A\subseteq \T$ on the
vector bundle $V\restr{X^{A}}$.  This makes it easy to study the effect on our
constructions of varying the representatives and for that matter the
subgroup $A$.  The second is a systematic use of the geometry of the
affine Weyl group of $\spin (2d)$, and its associated theory of theta functions
\cite{Looijenga:RootSystems}. 

More important than any single practical improvement was the
conceptual progress in understanding the relationship between the
sigma orientation and Looijenga's work on theta functions,
root systems, and elliptic curves.  In order to illustrate our
thinking, we state a conjecture; it is really a proposal for structure
which should be expected of equivariant elliptic cohomology, once a
rich enough theory has been found.  

Let $\Esym$ be the nonequivariant elliptic
spectrum associated to $C$ (see~\eqref{eq:19} for a construction), and
let 
\[
     \spfcohborel{X}{\Esym} \eqdef \spf E^{0} (\borel{X}).
\]
Then $X\mapsto \spfcohborel{X}{\Esym}$ is a covariant functor from
$\T$-spaces to formal schemes over $\spf E^{0} (B\T)\cong
\fmlgpof{C}$.
To state our conjecture, it is useful to suppose that $\T$-equivariant
elliptic cohomology is a covariant functor 
\[
   X \mapsto \ET{X}
\]
from $\T$-spaces to some category of
(super) ringed spaces over $C$, as was proposed by Grojnowski
\cite{Grojnowski:Ell} and Ginzburg-Kapranov-Vasserot \cite{GKV:Ell}
(and perhaps others).  For example, since Grojnowski's functor $\ETsym$
produces a $\Z/2$-graded commutative $\O{C}$-algebra, we may view
$(C,\ETS{X})$ as a ringed space $\ET{X}$ over $C$.    
In this notation the completion isomorphism~\eqref{eq:39} becomes
\[
     \spfcohborel{X}{\Esym} \cong (\ET{X} )^{\wedge}.
\]
Let $G = Spin (2d)$ with maximal
torus $T$ and Weyl group $W$, and let $\cochars=\hom (\T,T)$ be the
lattice of cocharacters.  The formal scheme $\cochars \otimes
\fmlgpof{C}$ carries a natural action of $W$, and the splitting
principle gives an isomorphism (see~\S\ref{sec-bg-e})
\[
    BG_{\Esym} \cong (\cochars\otimes \fmlgpof{C})/W.
\]
Suppose that $V$ is a $\T$-vector bundle over
$X$, and that $\borel{V}$ has structure group $G$.  The map 
\[
    X \xra{} BG
\]
classifying $V$ induces in $E$-cohomology a map
\begin{equation} \label{eq:38}
     \spfcohborel{X}{\Esym} \xra{} (\cochars\otimes \fmlgpof{C})/W.
\end{equation}
Looijenga constructs a holomorphic line bundle  $\mathcal{A}$ over
$(\cochars\otimes C)/W$ (\cite{Looijenga:RootSystems}; see also
\S\ref{sec:theta-functions}). The sigma function $\sigma$ defines a
global holomorphic section $\sigma_{d}$ of $\mathcal{A}$, whose zeroes
define an ideal sheaf $\I$ on $\VGW{C}$, such that $\sigma_{d}$ is a
\emph{trivialization} of $\mathcal{A}\otimes \I$.  Our conjecture is
the following; it is stated with a little more detail as Conjecture
\ref{t-co-princ-x-2}.

\begin{Conjecture} \label{t-co-princ-x}
The $\T$-equivariant elliptic cohomology $\ETsym$ associated
to the elliptic curve $C$ should associate to the $\T$-equivariant
spin bundle $V/X$ a map 
\[
    f:  \ET{X}  \rightarrow (\cochars\otimes C) /W, 
\]
which upon completion gives the map~\eqref{eq:38}.
Writing $\I (V) = f^{*} \I$ and $\mathcal{A} (V) = f^{*} \mathcal{A}$,
this map should have the following properties.
\begin{enumerate}
\item There is a canonical isomorphism 
\[
      \I (V) \cong \ETsym (V)
\]
of line bundles over $\ET{X}$.
\item Suppose that $V'$ is another $\T$-equivariant spin bundle.
If 
\begin{align*}
    c_{2} (\borel{V}) & = c_{2} (\borel{V'}), 
\end{align*}
then $\anomaly{V}\cong \anomaly{V'}$; indeed a $\T$-$\boeight$
structure on $V-V'$ determines a trivialization of $\anomaly{V}\otimes
\anomaly{V'}^{-1}$.
\end{enumerate}
\end{Conjecture}

In the notation of the Theorem, let $\sigma (V) = f^{*}
(\sigma_{d})$.  If $V- V'$ is a $\T$-$\boeight$-bundle, then $\sigma
(V)/ \sigma (V')$ is a trivialization of 
\[
   \frac{\anomaly{V}\otimes \I (V)}
        {\anomaly{V'} \otimes \I (V')}
   \cong 
   \frac{\ETS{V} }
        {\ETS{ V'}};
\]
this is the equivariant sigma orientation for $V-V'$.  From this point
of view, the sigma
function always determines a trivialization of the line bundle
$\anomaly{V}\otimes \I (V)$; when $c_{2} (V) = 0$, 
$\anomaly{V}$ is trivial and so this gives a Thom class.  The letter
$\mathcal{A}$ stands for ``anomaly''.  

Early versions of this paper were attempts to prove the Conjecture
\ref{t-co-princ-x} for Grojnowski's $\ETsym$, and then to deduce Theorem
\ref{t-th-sigma-intro} from it.  However, for reasons we explain in
\S\ref{sec:prin-bundles-ii}, we were only occasionally able to
convince even ourselves of those proofs.  Eventually, detailed
consideration of the consequences of the Conjecture led us to the
formulae in \S\ref{sec-useful-class} and so to concrete proofs.

In
\S\ref{sec:prin-bundles-ii} we do construct the map in the Conjecture
for a functor which captures the behavior of the stalks of
Grojnowski's functor.  The argument uses the same information we used
in \S\ref{sec-useful-class} and \S\ref{sec:thom-class} to produce the
equivariant Thom class.  The functor we study in
\S\ref{sec:prin-bundles-ii} was also inspired by Greenlees's rational
$\T$-equivariant elliptic spectra \cite{Greenlees:Ell} and by
Hopkins's work on characters in elliptic cohomology
\cite{Hopkins:Ell}.  Indeed we hope that Greenlees's rational
$\T$-equivariant elliptic spectra  will admit a proof
of the conjecture, and so give an account of the \emph{rational}
circle-equivariant sigma orientation.  

The rest of the paper proceeds as follows.  
In \S\ref{sec-notation} we summarize some of the notation which recurs
throughout the paper.  We discuss 
complex-orientable cohomology theories in general and ordinary and
elliptic cohomology theories in particular.  In 
\S\ref{sec-t-vbs} we state in a useful form some standard facts about
$\T$-equivariant principal $G$-bundles.

In~\S\ref{sec-bg-e} we interpret the analysis of \S\ref{sec-t-vbs} in
the presence of a (rational) complex-orientable cohomology theory.  We
begin \S\ref{sec:degr-four-char} with an interlude
(\S\ref{sec:degr-four-char-1}) on degree-four 
characteristic classes.  In \S\ref{sec:theta-functions} we
recall a result essentially due to \cite{Looijenga:RootSystems}, that
a degree-four characteristic class $\dfcc \in H^{4} (BG;\Z)$ gives rise
to a $W$-equivariant line bundle $\Loo (\dfcc )$ over
$(\cochars\otimes \CQ)$.  We define a \emph{theta function of level
$\dfcc$} for $G$ to be a $W$-invariant holomorphic section of the line
bundle 
$\Loo (\dfcc)$;   by the splitting principle, the Taylor series
expansion of such a theta function defines a characteristic class of
principal $G$-bundles.  The sigma function provides the most important
examples for us, and so we discuss it in 
\S\ref{sec:sigma-function-basic}.  

Section \ref{sec-useful-class} is the
heart of the paper.  In it we use the results of 
\S\ref{sec-t-vbs}---\ref{sec:sigma-function-basic} to construct some
holomorphic characteristic classes for $\T$-equivariant principal
$G$-bundles which are the building blocks of the Thom classes
in~\S\ref{sec:thom-class} and~\S\ref{sec:sigma-orientation}.

In~\S\ref{sec-equiv-ell} we recall the construction of Grojnowski's
analytic $\T$-equivariant elliptic cohomology associated to
a lattice $\Lambda\subset \C.$
In~\S\ref{sec:ellipt-cohom-thom-spaces}
we review the equivariant elliptic cohomology of Thom complexes  
\cite{Rosu:Rigidity,AndoBasterra:WGEEC}, recalling what is involved in
constructing a global section of $\ETS{V}$, where $V$ is (virtual)
$\T$-oriented vector bundle.

In~\S\ref{sec:sigma-orientation} we construct the equivariant sigma
orientation, proving Theorem \ref{t-th-sigma-intro}.
In~\S\ref{sec:thom-class} we prove the following related result.
Let $G$ be a spinor group, and let $G'$ be a simple and simply connected
compact Lie group.  Let $V$ be a $\T$-equivariant $G$-bundle 
over a finite $\T$-CW complex $X$, and let $V'$ be a $\T$-equivariant
$G'$-bundle.  Let 
\[
\Sigma (\borel{V}) \in  E (\borel{X}^{\borel{V}})
\]
be the Thom class given by the Weierstrass sigma function (see
Definition \ref{def-sg-or}).  Suppose 
that $\dfcc'$ is a degree-four characteristic classes for 
$G'$, with the property 
that 
\[
   c_{2} (\borel{V}) = \dfcc' (\borel{V'}).
\]
Suppose that $\theta'$ is a theta function for $G'$ of level
$\dfcc'$.  

\begin{Theorem}[\ref{t-th-thom-class}] \label{t-th-thom-class-intro}
A $\T$-orientation $\epsilon$ on $V$ determines a canonical global section
$\gamma=\gamma (V,V',\epsilon)$ of $\ETS{V}^{-1}$, whose 
value in $\stalk{\ETS{V}}^{-1}$ is $\theta' (\borel{V'})\Sigma
(\borel{V})^{-1}$.  The formation of $\gamma$ is natural in the sense
if 
\begin{align*}
     g/f: W/Y & \rightarrow V/X \\
     g'/f: W'/Y & \rightarrow V'/X
\end{align*}
are pull-backs of vector bundles, then under the isomorphism 
\begin{align*}
    \ETS{g/f}: \ETS{f}^{*} \ETS{V}  & \cong \ETS{W},
\intertext{we have}
    \ETS{f}^{*} \gamma (V,V',\epsilon) & = \gamma (W,W',g^{*}\epsilon).
\end{align*}
\end{Theorem}

In particular, suppose that $V'$ is an equivariant  $Spin (2d')$-vector bundle
over $X$, and $V$ is an equivariant $Spin (2d)$ vector bundle.   Suppose that
$\theta'$ is the character of a representation of $LSpin (2d')$ of
level $k$: then it is a theta function of level $k c_{2}$ for $Spin
(2d')$.  If
\[
    c_{2} (\borel{V}) = k c_{2} (\borel{V'}),
\]
then Theorem~\ref{t-th-thom-class-intro} gives a global section
$\gamma$ of
$\ETS{V}^{-1}$.  If $X$ is a manifold and $V$ is the tangent
bundle of $X$, then the Pontrjagin-Thom construction for the map
$\pi: X\to \point$ gives a  map 
\[
   (\ETS{\pi})_{*}  \ETS{V}^{-1} \rightarrow \ETS{* } = \O{C}
\]
of $\O{C}$-modules which takes $\gamma$ to the equivariant Witten
genus of $V$ twisted by the characteristic class $\theta'
(\borel{V'})$.  Since the global sections of $\O{C}$ are the
constants, we have the following result of Kefeng Liu \cite{Liu-1}.

\begin{Corollary}
Under these conditions, the equivariant Witten genus of $X$ twisted by
$\theta (\borel{V'})$ is constant.
\end{Corollary}

\begin{Remark} 
Liu states a condition on $p_{1}$ instead of $c_{2}$.
\end{Remark}

If \S\ref{sec-useful-class} is the heart of the paper, then
\S\ref{sec:prin-bundles-ii} is the soul.  There we discuss Conjecture
\ref{t-co-princ-x}, the study of which led to the results we report
here.  We give a refinement of the conjecture (Conjecture
\ref{t-co-princ-x-2}), and we explain how the arguments in this paper
support it.  We show how to construct the map of Conjecture
\ref{t-co-princ-x} for a functor which captures the behavior of the
stalks of Grojnowski's functor; the construction is
essentially a ``transfer formula'' in the sense of
\cite{BottTaubes:Rig}.   We show that the nonequivariant version of
the conjecture is true, and sheds light on the nonequivariant sigma
orientation.  Because the characters of representations of the loop group
$\Loops{G}$ are sections of the line bundle $\mathcal{A}$, it 
illuminates the relationship between the sigma orientation,
equivariant elliptic cohomology, and representations of loop groups
\cite{Brylinski:Ell,Ando:EllLG}. 

\subsection{Acknowledgments}  My first debt is to Maria Basterra.  Our
work on \cite{AndoBasterra:WGEEC} led directly to the results
in this paper, and I have very much enjoyed and benefitted from our
collaboration.  Maria was to have been an author of this paper as 
well, and I have reluctantly accepted her request to withdraw her name
from it.  I thank Alejandro Adem for inviting me to
visit Madison; it was during the eight-hour round-trip drive that the
the function $F$ \eqref{eq:23} was discovered, along with its
properties as described in \S\ref{sec-useful-class}.  I thank John Greenlees,
Haynes Miller, Jack Morava, and Charles Rezk for useful
conversations.  I thank Haynes Miller, Jack Morava, and Amnon Neeman
for encouraging me to record the ideas in
\S\ref{sec:prin-bundles-ii};  I hope they do not regret the result.
My work on elliptic cohomology has been profoundly influenced by Mike
Hopkins, and I am grateful to him for his work on the subject and his
generosity to me in particular. 

The author is supported by NSF grant DMS--0071482.  This paper was
completed during a visit to the Newton Institute for Mathematical
Sciences, during which time the author was supported by the Center for
Advanced Study of the University of Illinois at Urbana-Champaign.

\section{Notation}\label{sec-notation}

\subsection{Abelian groups}

Let $\cat{C}$ be a category with finite products.  The category $A \cat{C}$ of
abelian groups in $\cat{C}$ is an additive category.  In fact $A
\cat{C}$  is tensored over the category of finitely
generated free abelian groups. 
That is, a finitely generated free abelian group $F$ and an abelian
group $X$ of $\cat{C}$ determine (naturally in $F$ and $X$) an object
$F\otimes X$ of $A \cat{C}$, with a natural isomorphism
\[
 A\cat{C} [F\otimes X, Y]\cong \CategoryOf{abelian groups}[F,A\cat{C}[X,Y]].
\]
If $A$ is an abelian group written additively, and $M$ is an abelian
group written multiplicatively, then we write $m^{a}$ for the element
$a\otimes m$ of $A\otimes M$.  Similarly, if $M'$ is an
abelian group, $u\in M'$, and $m\in \hom (M',M)$ then we may write $u^{m}$
for $m (u)$.

We write $\Ga$ for the additive group, and $\Gm$ for the
multiplicative group.

\subsection{Lie groups and the group $\VG{X}$}

In general, the letter $G$ will stand for a compact Lie group with maximal
torus $T$ with Weyl group $W$.  We define 
\begin{align*}
\cochars&\eqdef\hom[\T,T] \\
\chars  &\eqdef\hom[T,\T]
\end{align*}
to be the lattices of cocharacters and characters.   We write 
\[
    \Ad: G\to \Aut (G)
\]
for the action of $G$ on itself by \emph{conjugation}: 
\[
    \Ad_{g} h = g h g^{-1}.
\]
If $X$ is an abelian group in any category, then the tensor product
$\VG{X}$ carries an action of the Weyl group $W$.  If $\rank$ is the
rank of $G$, then $\VG{X}$ is isomorphic to $X^{r}$.

\subsection{Elliptic curves}
\label{sec:elliptic-curves}

Fix  $\tau$ in the complex upper half plane, and let $\Lambda$ be the
lattice 
\[
 \Lambda = 2 \pi i \Z + 2\pi i \tau \Z.
\]
The complex numbers appear in this paper both as a ring and as an
analytic variety.  To avoid confusion we write $\affan$ for the
complex numbers regarded as an analytic variety.  Similarly, we write
$\Gaan$ for $\affan$ with its additive structure of abelian
topological group.  We write $z$ for the 
standard coordinate on $\Ga$ and also on $\aff$, $\affan$, $\Gaan$,
etc.  We write $\Gman$ for 
$(\affan)^{\times}$ with its multiplicative group structure.  
We set
\begin{equation} \label{eq:43}
\begin{split}
  u^{r} & = e^{r z} \\
  q^{r} & = e^{2\pi i r \tau} 
\end{split}
\end{equation}
for $r\in \Q$, and 
we let $\CQ$ be the elliptic curve
\[
  \CQ = \C / \Lambda = \Gaan/\Lambda \cong \C^{\times} / q^{\Z} =
\Gman /q^{\Z}.
\]
We write $\pr$ for the covering map 
\[
    \affan \xra{\pr} \CQ.
\]
If $V$ is an open set in a  complex analytic variety, then we write
$\O{V}$ for the sheaf of holomorphic functions on $V$.

If $A$ is an abelian
topological group and $a\in A$, when we write $\trans_{a}$ for the
translation map; and if $V\subset A$ is an open set, then we write
\[
    V-g \eqdef \trans_{-g} (V).
\]

\begin{Definition} \label{def-small}
An open set  $U$ of 
$C$ is \emph{small} if it is connected and $\pr^{-1}U$ is a union of
connected components $V$ with the property that 
\[
      \pr\restr{V}: V\to U
\]
is an isomorphism.
\end{Definition}  

If $U$ is small and $V$ is a component of $\pr^{-1}U$, then 
the covering map induces an isomorphism 
\[
    \O{U}\cong \O{V}.
\]
In particular, if $U$ contains the origin of $C$, then there
is a unique component $V$ of $\pr^{-1}U$ containing $0$.  This
determines a $\C[z]$-algebra structure on $\O{U}$, and a
$\C[z,z^{-1}]$ structure on  $\O{U}\restr{U\backslash 0}.$

\subsection{Ringed spaces}

Grojnowski's $\T$-equivariant elliptic cohomology is a contravariant
functor 
\begin{equation} \label{eq:25}
 \ETsym: \CategoryOf{finite $\T$-CW complexes} 
\rightarrow 
 \CategoryOf{$\Z/2$-graded $\O{C}$-algebras}
\end{equation}
(see \S\ref{sec-equiv-ell} and
\cite{Grojnowski:Ell,Rosu:Rigidity,AndoBasterra:WGEEC}).  At roughly
the same time as \cite{Grojnowski:Ell}, Ginzburg-Kapranov-Vasserot
\cite{GKV:Ell} proposed that the $\T$-equivariant elliptic cohomology
associated to an elliptic curve $C$  should be a \emph{covariant} functor 
\[
   E_{\gkv}: \CategoryOf{$\T$-spaces}\rightarrow \CategoryOf{schemes}/C,
\]
although they gave no independent construction of such a functor.  

These are
meant to be related by the formula 
\begin{equation} \label{eq:26}
    \ETS{X} = f_{*}\O{E_{\gkv} (X)},
\end{equation}
where 
\[
   f : E_{\gkv} (X) \rightarrow C
\]
is the structural map.  Grojnowksi's functor can not quite be
of the 
form~\eqref{eq:26}, since in~\eqref{eq:25} $\O{C}$ is the sheaf of 
\emph{holomorphic} functions on the \emph{analytic} space
$C=\Gaan /\Lambda$.    However, it does give a covariant functor 
\[
\CategoryOf{finite $\T$-spaces}\rightarrow \CategoryOf{ringed spaces}/C.    
\]
Precisely, we have the following.

\begin{Definition}\label{def-ringed-space} 
By a \emph{(super, or $\Z/2$-graded) ringed space} we shall mean a
pair $(X,\O{X})$ consisting of a space $X$ and a sheaf $\O{X}$ of 
$\Z/2$-graded rings on $X$.  A \emph{map of ringed spaces}
\[
f= (f_{1},f_{2}) : (X,\O{X}) \rightarrow (Y,\O{Y})
\]
consists of a map of spaces $f_{1}: X\to Y$ and a map of sheaves of
$\Z/2$-graded commutative algebras over $Y$
\[
        f_{2}: \O{Y} \rightarrow ( f_{1})_{*} \O{X}.
\]
The resulting category of ringed spaces will be denoted
$\ringedspaces$.
If $\mathcal{X} = (X,\O{X})$ is a
ringed space and $U$ is an open set of $X$, then we may write
$\mathcal{X} (U)$ in place of  $\O{X}(U)$.
\end{Definition}

If $X$ is a finite $\T$-CW complex, then 
\[
\ET{X} \eqdef  (C,\ETS{X})
\]
is a ringed space.  If 
\[
  f : Y\to X
\]
is a map of finite $\T$-CW complexes then we define 
\[
   \ET{f}: \ET{Y} \to \ET{X}
\]
to be the identity on $C$ and 
\[
\ETS{f}: \ETS{X} \rightarrow \ETS{Y}
\]
on the structure sheaves.  In this we obtain a covariant functor
\[
\spaceof{\ETsym} : \CategoryOf{finite $\T$-CW complexes}
\rightarrow \ringedspaces/C.
\]
We have found this point of view to be extremely helpful, and so we
have adopted it in writing this paper.


\subsection{Ordinary cohomology}

If $R$ is a commutative ring, let $HR$ denote ordinary cohomology
with coefficients in $R$.  If $X$ is a space, then we write
\[
     \spaceof{HR} (X) \eqdef \spec (HR^{\text{even}} (X))
\]
for the scheme over $\spec R$ associated to the cohomology of $X$.  We may
also view $\spaceof{HR} (X)$ as the ringed space with underlying space $\spec
(HR^{\text{even}} (X))$  and structure sheaf associated to 
\[
     HR^{\text{even}} (X)\oplus HR^{\text{odd}} (X).
\]
We shall write $H$ for $H\C$, cohomology with complex coefficients.

\subsection{Equivariant cohomology}

If $X$ is a space with a circle action, 
$\spaceof{HR}(\borel{X})$
is a scheme over $\spaceof{HR} ( B\T) $, which we denote $\spaceof{HR} (X).$
We choose a generator of the character group of $\T$, and write $z$
for the resulting generator of $\HZsym ^{2} (B\T)$; this gives an
isomorphism 
\begin{equation}\label{eq:12}
     \spaceof{HR} ( B\T)  \cong (\Ga)_{R}
\end{equation}
of group schemes over $\spec R.$  We shall use \eqref{eq:12} to view
$\spaceof{HR} (\borel{X})$ as a scheme over $\aff_{R}$. 

We recall \cite{Quillen:secr} that 
equivariant cohomology satisfies a localization theorem.

\begin{Theorem} \label{t-th-localization}
If $X$ has the homotopy type of a finite $\T$-CW complex (e.g. if $X$
is a compact $\T$-manifold), then the natural map 
\[
     \THC{X^{\T}} \xra{} \THC{X}
\] 
induces an isomorphism over $\spec \C[z,z^{-1}] \subset \HC{B\T}$.
\qed
\end{Theorem}

\subsubsection{Holomorphic cohomology}

Let $\affan$ be the analytic complex plane, so $\O{\affan}$ is the
sheaf of holomorphic functions on $\C$.  Because of the natural maps 
\[
   \affan \rightarrow \aff_{\C} \rightarrow \aff \cong \HZ{B\T}
\] 
we may view $z$ as a function on $\affan$.  Given a $\T$-space $X$ we
define the \emph{holomorphic} cohomology of $X$ to be the 
sheaf of super $\O{\affan}$-algebras given by 
\[
    \HHol{X}{U} \eqdef H \C (\borel{X}) \tensor{\C[z]}
\O{\affan} (U). 
\]
We view $\HHsym (X)$ as the structure sheaf of a ringed
space~\eqref{def-ringed-space} $\TH{X}$ over $\affan$, namely the the
pull-back in the diagram 
\[
\begin{CD}
\TH{X} @>>> \THC{X} @>>> \THZ{X} \\
@VVV       @VVV           @VVV \\
\affan @>>> \aff_{\C} @>>> \aff_{\Z}.
\end{CD}
\]

\subsubsection{The stalk of holomorphic cohomology}

Let 
\[
\affanstalk = \spec (\O{\affan,0})
\]
be the local scheme associated to the stalk of $\O{\affan}$ at the
origin.  The stalk of $\HHsym (X)$ at the origin is 
\[
   \HHstalk{X}  \iso 
   \HCsym (\borel{X})\tensor{\C[z]} (\O{\affan,0}).
\]
We write $\THstalk{X}$ for the resulting scheme over $\affanstalk$.

\subsubsection{Periodic Borel cohomology}

Let $\HPCsym$ denote \emph{periodic} ordinary cohomology with complex
coefficients: that is, 
\begin{equation} \label{eq:22}
     \HPCsym = \bigvee_{k\in \Z} \Sigma^{2k} \HCsym
\end{equation}
so 
\[
     \pi_{*} \HPCsym = \C[v,v^{-1}]
\]
with $v \in \pi_{2} \HPCsym$.    Then 
\[
     \HPCsym^{0} (B\T) \iso \C\psb{z}
\]
so $\spf \HPCsym^{0} (B\T) = \affanfml = (\aff_{\C})^{\wedge}_{0}$,
and $\HPCsym^{0} (X_{\T})$ is the ring of formal functions on the
pull-back $\THPC{X}$ 
in the diagram of formal schemes
\[
\begin{CD}
\THPC{X} @>>> \THstalk{X}\\
@VVV @VVV \\
\affanfml @>>> \affanstalk.
\end{CD}    
\]
Summarizing, we have a collection of forms of ordinary
cohomology 
\begin{equation}\small \label{eq:2}
\begin{CD}
\THPC{X} @>>> \THstalk{X} @>>> \TH{X} @>>> \THC{X} @>>> \THZ{X} \\
@VVV          @VVV              @VVV        @VVV        @VVV \\
(\aff_{\C})^{\wedge}_{0} @>>>
           \affanstalk    @>>> \affan @>>>  \aff_{\C}  @>>> \aff_{\Z}.
\end{CD}
\end{equation}

\subsection{Generalized cohomology}

\subsubsection{Even periodic ring spectra}

A ring spectrum $\Esym$ will be called ``even periodic'' if
$\pi_{\text{odd}}\Esym=0$  and $\pi_{2}\Esym$ contains a unit of
$\pi_{*}\Esym$.  

If $\Esym$ is an even periodic ring spectrum, and if $X$ is a space, then
we shall write  $\Esym^{*} X$ for the unreduced cohomology of $X$.  
As in \cite{AHS:ESWGTC}, we write $\E{X}$ for the formal scheme 
\[
   \E{X} = \colim_{F\subset X} \spec \Esym^{0}F
\]
over $S_{\Esym} = \spec \Esym^{0} (\point)$; the colimit is 
over the compact subsets of $X$.

An even periodic ring spectrum is always complex-orientable.  In particular
\[
\GpOf{\Esym} \eqdef \E{B\T}
\]
is a (commutative, one-dimensional) formal group over $S_{\Esym}$.
For example, let $\HPZsym$ denote periodic ordinary cohomology with
integer coefficients.  Then $\GpOf{\HPZsym} = \Gah$.

\subsubsection{Borel cohomology}

If $\Esym$ is an even periodic ring spectrum, and 
$X$ is a space with a circle action, then the projection
\[
\borel{X}\rightarrow B\T   
\]
induces a map 
\[
\spfcohborel{X}{\Esym} \rightarrow \GpOf{\Esym},
\]
making
\[
\spfcohborel{X}{\Esym} = \spf \Esym^{0} (\borel{X})
\]
into a formal scheme over the formal group of $\Esym$.

\subsubsection{Elliptic spectra}

We recall \cite{AHS:ESWGTC} that an \emph{elliptic spectrum} is a
triple $(E,C,t)$ consisting of  
\begin{enumerate}
\item an even periodic ring spectrum $E$,
\item a (generalized) elliptic curve $C$ over $S_{E}$, and
\item an isomorphism of formal groups 
\[
    t: \fmlgpof{C}\cong \GpOf{E}.
\]
\end{enumerate}

\subsubsection{Rational elliptic spectra}

Let 
\[
\xymatrix{
{C}
   \ar[r] 
&
{\spec R}
 \ar@/_1pc/[l]_{0}
}
\]
be an elliptic curve with identity $0$ over a $\Q$-algebra $R$, and let
$\oneforms = 0^{*}\Omega^{1}_{C/R}$.  If $\omega$ is a trivialization
of $\oneforms$, then there is a canonical
isomorphism of formal groups over $R$
\[
      \fmlgpof{C}\xrightarrow[\cong]{\log_{\omega}} \Gah
\]
with the property that $0^{*}\log_{\omega}^{*}dz = 0^{*}dz$.

Let $\Gamma^{\times} (\oneforms)$ be the functor from rings to
sets which given by 
\[
   \Gamma^{\times} (\oneforms) (T) = \{(j,\omega) \suchthat j: \spec T\to
\spec R, \omega \text{ a trivializaton of }j^{*}\oneforms \}.
\]
If $\omega$ is a trivialization of $\oneforms$, then the
trivializations of $j^{*}\oneforms$ are of the form $u j^{*}\omega$ for
$u\in T^{\times}$, so 
\[
    \Gamma^{\times} (\oneforms) \iso \spec R[u,u^{-1}].
\]
Let $S=\O{}(\Gamma^{\times} (\oneforms))$. 

The logarithm then gives a canonical isomorphism of formal groups
\[
    \fmlgpof{C}_{S} \xra{\log_{\fmlgpof{C}}} (\Gah)_{S}
\]
by the formula 
\[
      (c,\omega) \mapsto (\log_{\omega}c, \omega).
\]
If $HS[v,v^{-1}]$ is the even periodic ring spectrum such that 
\[
    (HS[v,v^{-1}])^{*} X \eqdef H^{*} (X;S[v,v^{-1}]),
\]
where we take $v$ to have degree $2$ and $u\in S$ to have degree zero,
then 
\[
    \GpOf{HS[v,v^{-1}]} = (\Gah)_{S},
\]
and we have an elliptic spectrum
\[
     (HS[v,v^{-1}],C,\log_{\fmlgpof{C}}).
\]
Alternatively, given over $R$ a trivialization $\omega$ of $\oneforms$ we have 
the elliptic spectrum $(HR[v,v^{-1}],C,\log_{\omega})$.

\subsubsection{Complex elliptic spectra}
 
Recall~\eqref{eq:22} that $\HPCsym$ denotes periodic ordinary
cohomology with complex coefficients.

The projection $\pr: \Gaan \to C$ induces an isomorphism of formal
groups 
\[
\fmlgpof{\pr}: \Gah \to \fmlgpof{C}.
\]
There is a unique cotangent vector $\omega$ such that 
\[
   \pr^{*}\omega = 0^{*}dz.
\]
We have 
\[
    \log_{\omega} = (\fmlgpof{\pr})^{-1},
\]
and so an elliptic spectrum 
\begin{equation} \label{eq:19}
(\HPCsym,
C,\fmlgpof{\pr}^{-1}) = (\HPCsym,C,\log_{\omega}).
\end{equation}

\section{Principal bundles with an action of the circle}
\label{sec-t-vbs}

\subsection{$A$-bundles over trivial $A$-spaces}

Let $A$ be a closed subgroup of the circle $\T$.  Suppose that $G$ is
a \emph{connected} compact Lie group with maximal torus $T$, and let
$\pi: \PB\to Y$ be a principal $G$-bundle 
over a connected space $Y$.  Suppose that $A$ acts on $\PB/Y$,
fixing $Y$. The group of automorphisms of $\PB/Y$ is the group of
sections $\Gamma (( \PB\times_{G} G^{\Ad})/Y)$, and an action of $A$ on
$\PB/Y$ is equivalent to a section $a \in \Gamma ((\PB\times_{G}
\hom[A,G^{\Ad}])/Y)$ (The notation $G^{\Ad}$ refers to $G$ with $G$
acting on it by conjugation).

Since $G$ is connected and $A$ is (topologically) cyclic, every
$G$-orbit in $\hom[A,G]$ intersects   $\hom[A,T]$ nontrivially: that
is, the map 
\[
   \hom[A,T] \rightarrow \hom[A,G^{\Ad}] \rightarrow \hom[A,G^{\Ad}]/G
\]
is surjective. 

In particular $\hom[A,G^{\Ad}]/G$ is discrete.  Its points label the
connected components of  $\PB\times_{G} \hom[A,G^{\Ad}]$ via the
surjective map 
\[
      \PB\times_{G} \hom[A,G^{\Ad}] \rightarrow \hom[A,G]/G.
\]
Since $Y$ is connected, we may choose a homomorphism 
$m\in \hom[A,T]$ such that, for all $x\in Y$, there is a $p\in \PB$
such that 
\[
    a (x) = [p,m];
\]
the square brackets indicate the class in the Borel construction of the 
element $(p,m)\in \PB\times\hom[A,G]$.  The choice of $m$
determines  $p$ only up to the centralizer $Z (m)$ in $G$ of the
homomorphism $m$.

\begin{Definition} \label{def-reduction}
A \emph{reduction} of the action of $A$ on $\PB$ is a homomorphism 
\[
    m: A\to T
\]
such that, for all $x\in Y$, there is a $p\in \PB$ such that 
\[
 a (x) = [p,m].
\]
\end{Definition}
 
The terminology is justified by the following observation. 
Let 
\[
   \PB (m) = \{p\in \PB| [p,m] = a (\pi (p)) \}.
\]
Then $\pi\restr{\PB (m)}:\PB (m)\to Y$ is a principal $Z
(m)$ bundle over $Y$, 
the reduction of the structure group of $\PB$ to $Z (m)$.  In other
words, we have given a factorization 
\begin{equation}\label{eq-red-structure-gp}
\xymatrix{
&
{BZ (m)}
  \ar[d]\\
{Y}
 \ar[r]_{\PB}
 \ar[ur]^{\PB (m)}
&
BG.
}
\end{equation}
Let 
\[
W (m) = \{w \in W| w m = m\}
\]
be the stabilizer of $m$.  One sees that $T$ is a maximal
torus of $Z (m)$, with Weyl group $W (m)$.   

\begin{Example}\label{ex-U-n-reduction}
Suppose that $G=U (n)$ is the unitary group with its maximal torus
$T=\Delta (z_{1},\dots ,z_{n})$ of diagonal matrices, and suppose that
$A=\T$.  Every homomorphism
\[
     m: \T\to T
\]
is conjugate in $U (n)$ to one of the form 
\[
    m (z) = \Delta (z^{m_{1}},\dots ,z^{m_{1}},z^{m_{2}},\dots
,z^{m_{2}},\dots ,z^{m_{r}},\dots ,z^{m_{r}}), 
\]
where $m_{i}\in \Z$.  Let $d_{j}$ be the multiplicity of $m_{j}$; then
the centralizer is the block-diagonal matrix 
\[
   Z (m) = \left[
\begin{array}{ccc}
U (d_{1}) & & \\
          &\dots &\\
          & & U (d_{r})
\end{array} \right],
\]
with Weyl group
\[
   W (m) = \Sigma_{d_{1}}\times\dots \times \Sigma_{d_{r}} \subset
\Sigma_{n} = W.
\]
If $V$ is a $\T$-equivariant complex vector bundle over a connected trivial
$\T$-space $Y$, and if this  $m$ is a reduction of the action of $\T$
on the principal bundle of $V$, then the reduction of the structure
group to $Z (m)$ corresponds to the decomposition of $V$ as the direct sum
\[
    V \cong V (m_{1})\oplus \dots \oplus V (m_{r}),
\]
where $\T$ acts on the fiber of $V (m_{j})$ by the character
$z^{m_{j}}$.  
\end{Example}

The composition 
\[
 g (m):  A \times Z (m)
        \xrightarrow{m \times Z (m)}
        T\times Z (m) \rightarrow Z (m)
\]
is a group homomorphism, and so we have a map 
\[
   BA  \times BZ (m) \xra{Bg (m)} BZ (m).
\]
The following Lemma will be used directly to prove Lemma
\ref{t-le-equiv-char-class-with-I-and-qfm}, a calculation of
degree-four characteristic classes.  Moreover, the algebro-geometric form
of this diagram after applying a complex-orientable cohomology theory
(see Lemma \ref{t-le-equiv-char-class-restr-VC}) captures the
essential point of the ``transfer argument'' of
\cite{BottTaubes:Rig}, as we explain in Remark \ref{rem-1}.

\begin{Lemma}  \label{t-le-g-m}
\begin{enumerate}
\item 
The diagram 
\[
\xymatrix{
&
{BT\times BT}
 \ar[dr]
\\
{BA \times BT}
 \ar[ur]^{Bm\times BT}
 \ar[d] 
&&
{BT}
 \ar[d] \\
{BA\times BZ (m)} 
  \ar[rr]^{Bg (m)}
&&
{BZ (m)}
}
\]
commutes. 
\item 
The map $Bg (m)$ classifies the principal $Z (m)$-bundle $EZ (m)_{A} =
EA\times_{A} EZ (m)$ over $BA\times BZ (m)$.
\end{enumerate}
\end{Lemma}

\begin{proof}
The first part is easy.  
For the second part, it suffices to construct 
a map of principal $Z (m)$-bundles over $BA\times BZ (m)$.
The map 
\[
   Eg (m) : EA \times EZ (m) \rightarrow EZ (m)
\]
factors through $EA\times_{A} EZ (m)$, and gives the desired map.
(Note that the map $Eg (m)$ is obtained by constructing  $EA$ and $EZ
(m)$ functorially as spaces with actions on the
same side, say the left.  In forming $EA\times_{A} EZ (m)$, one
makes $A$ act on the right of $EA$ by the inverse.)
\end{proof}

If $m': A\to T$ is another reduction of the action of $A$ on $Q$, then
$m$ and $m'$ differ by conjugation in $G$, and we have for some $g$ in
$G$ a commutative diagram 
\begin{equation} \label{eq:1}
\xymatrix{
{Y}
\ar[rr]^{\PB (m')}
\ar[dd]_{\PB (m)}
\ar[ddrr]_>>>>>>>>{\PB}
& & 
{BZ (m')} 
 \ar[dd]
\\
 & 
 & 
\\
{BZ (m)} \ar[rr] \ar '[ur]^{c_{g}} [uurr]
& & 
{BG.}
}
\end{equation}

\subsection{The case of a connected centralizer}

\begin{Lemma} \label{t-le-g-in-NT-if-Z-conn}
If the centralizer $Z (m)$ is
connected, then the element $g$ in the diagram \eqref{eq:1} may be
taken to be in the normalizer $N_{G}T$ of $T$ in $G$.
\end{Lemma}

\proof
Let $g\in G$ be such that 
\[
    \Ad_{g} m = m'
\]
Then $m' (A)\subset \Ad_{g} (T)\cap T$, so both $T$ and
$\Ad_{g} (T)$ are maximal tori in $Z (m')$.  Since $Z (m')$ is
connected, there is an element $h\in Z (m')$ such that 
\[
    \Ad_{h}\Ad_{g} (T) = T,
\]
so $hg\in N_{G} T$.  Since $h\in Z (m')$, we have 
$$
m' = \Ad_{h}m' = \Ad_{h} \Ad_{g} m.\eqno{\qed}
$$

Example~\ref{ex-U-n-reduction} shows that $Z (m)$ is connected if $G$
is a unitary group.  Bott and Samelson have shown
(\cite{BS:atmss}; see Proposition 10.2 of \cite{BottTaubes:Rig})
that $Z (m)$ is connected if $G$ is a spinor group.  Indeed they prove
the following.

\begin{Lemma}
\label{t-le-Z-connected} 
If $G$ is simple and simply connected then the centralizer $Z (m)$ is
connected. \qed
\end{Lemma}

\subsection{Nested fixed-point sets}

Now suppose that $B\supseteq A$ is a larger closed subgroup of $\T$
(primarily, we shall be interested in the case that $B=\T$), that
$\PB/Y$ is a $B$-equivariant principal $G$-bundle, and that $A$ acts
trivially on $Y$.

\begin{Lemma} \label{t-le-borel-of-reduced}
If 
\[
m: A\to T
\]
is a reduction of the action of $A$ on $\PB$, then the action of $B$ on
$\PB$ induces on $\PB (m)$ the structure of a $B$-equivariant 
$Z (m)$-bundle over $Y$.  
Moreover the Borel
construction $\PB (m)_{B}$  is
a principal $Z (m)$-bundle over $Y_{B}$.   
\end{Lemma}

\begin{proof}
This follows from the fact that $B$ is abelian.
\end{proof}

If $F\subseteq Y^{B}$ is a connected component of the subspace of $Y$ fixed by
the action of $B$, then we may choose a reduction
\[
    m_{F}: B\to T
\]
of the action of $B$ on $\PB\restr{F}$.  

\begin{Lemma} \label{t-le-m-Y-restr-to-Fcheck}
The restriction 
\[
m_{F}\restr{A}: A\to T
\]
is a reduction of the action of $A$ on $Y$.
\end{Lemma}

\proof
The action of $A$ on $\PB/Y$ is a section $a$ of $(\PB\times_{G}
\hom[A,G^{\Ad}])/Y$.  The restriction $a\restr{F}$ records the action of $A$ on
$\PB\restr{F}$.  The action of $B$ on $\PB\restr{F}$ is a section $b$ of 
$(\PB\restr{F} \times_{G} \hom[B,G^{\Ad}])/F$, and with the obvious notation we
have 
$$
     b\restr{A} = a\restr{F}.\eqno{\qed}
$$

\section{Rational cohomology of principal bundles with\nl compact connected
structure group}
\label{sec-bg-e}

Let $\Esym$ be an \emph{rational} even periodic ring spectrum.
Let $G$ be a connected compact Lie group.  Let $T$ be a maximal torus
of $G$, with Weyl group $W$.  The natural isomorphism
\[
         \cochars \otimes \T \to T
\]
induces a $W$-equivariant isomorphism
\[
         \VG{\GpOf{\Esym}} = \VG{\E{B\T}}  \cong \E{BT}
\]
of formal groups over $S_{\Esym}$.  Moreover \cite{Borel:tlgcc} the
natural map  
\[
         \E{BT} / W \to \E{BG}
\]
is an isomorphism.  We shall repeatedly use the
resulting isomorphism 
\[
         \E{BG} \cong \VGW{\GpOf{\Esym}}.
\]
For example, a principal $G$-bundle $\PB$ over $X$ is classified by a
map 
\[
     X\xra{\PB} BG
\]
whose effect in $\Esym$-theory 
\begin{equation}\label{eq:47}
    \E{X} \xra{\E{\PB}} \VGW{\GpOf{\Esym}}
\end{equation}
is a $\E{X}$-valued point of $\VGW{\GpOf{\Esym}}$.

\subsection{Periodic cohomology of circle-equivariant principal bundles}

We interpret the analysis of \S\ref{sec-t-vbs} in $\Esym$-theory.
Suppose that $G$ is a connected compact Lie
group, and that $\PB$ is a $\T$-equivariant principal $G$-bundle over
a connected $\T$-space $Y$.  Suppose that a closed subgroup $A$ of the circle 
acts trivially on $Y$.  Suppose that 
\[
  m: A\to T
\]
is a reduction of the action of $A$ on $\PB/Y$ with \emph{connected}
centralizer  $Z (m)$.

Applying $\Esym$-cohomology to the diagram \eqref{eq-red-structure-gp} yields
\[
\xymatrix{
&
{\VG{\GpOf{\Esym}}/W (m)}
  \ar[d]\\
{\E{Y}}
 \ar[r]_-{\E{\PB}}
 \ar[ur]^-{\E{\PB (m)}}
&
{\VG{\GpOf{\Esym}}/W}
}
\]

The multiplication 
\[
    BT\times BT \rightarrow BT
\]
induces the addition map 
\[
   (\VG{\GpOf{\Esym}})\times  (\VG{\GpOf{\Esym}})
    \xrightarrow{+} \VG{\GpOf{\Esym}}
\]
whose restriction to 
\[
    (\VG{\GpOf{\Esym}})^{W (m)}\times  (\VG{\GpOf{\Esym}})
    \xrightarrow{+} \VG{\GpOf{\Esym}}
\]
factors to give a translation map 
\[
    (\VG{\GpOf{\Esym}})^{W (m)} \times  (\VG{\GpOf{\Esym}})/W (m)
    \xrightarrow{+} (\VG{\GpOf{\Esym}})/ W (m).
\]
In $\Esym$ cohomology, Lemma \ref{t-le-g-m} implies the following
result.  It is used in the case $A=\T$ to construct the commutative
diagram \eqref{eq:4} and so prove Lemma \ref{t-le-translation-on-F}.
It is also implies the commutativity of the diagram
\eqref{eq-trans-form-diagram}, which captures the essence of the
``transfer formula'' of Bott-Taubes \cite{BottTaubes:Rig}; see Remark
\ref{rem-1}.

\begin{Lemma} \label{t-le-equiv-char-class-restr-VC}
The diagram 
\[
\xymatrix{
{\E{BA} \times \E{BZ (m)}} 
 \ar[r]
&
{(\VG{\GpOf{\Esym}})^{W (m)} \times (\VG{\GpOf{\Esym}})/W (m)}
 \ar[d]^-{+}
\\
{\E{BA} \times \E{Y}}
 \ar[u]^-{\E{BA}\times \E{Q (m)}}
 \ar[r]^-{\E{\PB (m)_{A}}}
 \ar[d]
&
{(\VG{\GpOf{\Esym}})/W (m)} 
\\
{\spfcohborel{Y}{\Esym}}
\ar[ur]_-{\spfcohborel{\PB (m)}{\Esym}}
}
\]
commutes.  
\qed
\end{Lemma}

\subsection{Holomorphic characteristic classes}
\label{sec:holom-char-class}

Let $f$ be a $W$-invariant holomorphic function on $\VG{\Gaan}$.  By
the splitting principle, the Taylor expansion of $f$ at the origin
defines a class in $\HPCsym (BG)$.  Suppose that $\PB$ is a principal
$G$-bundle over $X$, with the property that $\borel{\PB}$ is a
principal $G$-bundle over $\borel{X}$.  Then we get a class
\[
      f (\borel{\PB})  \in \HPCsym (\borel{X}).
\]
The following result is due to Rosu.
\begin{Lemma} \label{t-le-hol-char-classes-I}
The class $f (\borel{\PB})$ is in fact an element of
$\HHol{X}{\affan}$.  Similarly, if $f\in \stalk{(\O{\VG{\Gaan}})}$, then
\[
     f (\borel{\PB }) \in \HHstalk{X}.
\]
\end{Lemma}

\begin{proof}
Proposition A.6 of \cite{Rosu:Rigidity} proves this result in the case
that $G=U (n)$ and 
\[
f (z) = \prod_{j} g (z_{j}),
\]
where $g\in \O{\affan,0}$ and 
$z = (z_{1},\dots ,z_{n}) \in \mathbb{A}_{\mathrm{an}}^{n} \cong
\VG{\Gaan}$.  The same argument works in the indicated generality.
\end{proof}

An important example of a holomorphic characteristic class is the
Euler class associated to a ``multiplicative analytic orientation''.  
A power series
\[
      f (z) = z + \text{ higher terms} \in \HPCsym (B\T) \cong \C\psb{z}
\]
satisfying 
\[
     f (-z) = - f (z)
\]
determines a multiplicative orientation  (map of homotopy commutative
ring spectra)
\[
    \phi: MSO \rightarrow \HPCsym,
\]
characterized by the property that if 
\[
V = L_{1} + \dots  + L_{d}
\]
is a sum of complex line bundles, 
then its Euler class is 
\[
  e_{\phi} (V) = \prod_{j} f (c_{1}L_{j}).
\]

If $V$ is  an oriented vector bundle, we write $\phi (V)\in \HPCsym
(V) = \HPCsym (X^{V})$ 
for the resulting Thom class.  It is multiplicative in the sense that   
\begin{equation}\label{eq-thom-iso-mult}
    \phi (V \oplus W) = \phi (V) \wedge \phi (W)
\end{equation}
under the isomorphism 
\[
     (X\times Y)^{V \oplus W} \cong X^{V}\wedge Y^{W}.
\]

\begin{Definition}\label{def-analytic-or} 
The orientation $\phi$ is \emph{analytic}
if $f$ is contained in the subring
$\O{\affan,0} \subset \C\psb{z}$ of germs of holomorphic functions
at $0$; equivalently, if there is a neighborhood $U$ of $0$ in $\C$ on
which the power series $f$ converges to a holomorphic function.  
\end{Definition}
Lemma~\ref{t-le-hol-char-classes-I} implies the following.

\begin{Corollary}{\rm\cite{Rosu:Rigidity}}\qua\label{t-co-analytic-euler-classes} 
If $\phi$ is analytic and $V$ is an oriented $\T$-vector bundle over a
compact $\T$-space $X$, then the Euler class $e_{\phi}$
associated to $\phi$ satisfies 
$$
   e_{\phi}(\borel{V}) \in \HHstalk{X}.\eqno{\qed}
$$
\end{Corollary}

If $\Phi$ denotes the standard Thom isomorphism, then 
\[
   \phi (\borel{V}) = \frac{e_{\phi} (\borel{V})}
                     {e_{\Phi} (\borel{V})} \Phi (\borel{V}),
\]
and the ratio of Euler classes is a unit in $\HHstalk{X}.$ 
Of course multiplication by $\Phi (\borel{V})$ induces an isomorphism 
\[
    \HCsym (\borel{X}) \cong \HCsym (\borel{V}),
\]
and so we have the following.

\begin{Corollary} \label{t-cor-analytic-thom}
There is a neighborhood $U$ of the origin in $\affan$ such
that 
\[
    \phi (\borel{V}) \in \HHol{V}{U},
\]
and such that multiplication by this class induces an isomorphism of
sheaves 
\[
   \HHsym (X)\restr{U}  \xrightarrow[\cong]{\phi}    
   \HHsym (V)\restr{U}.
\]
In other words, for every open set $U'\subseteq U$, multiplication by
$\phi (\borel{V})$ induces an isomorphism 
$$
   \HHol{X}{U'}  \xrightarrow[\cong]{\phi}    
   \HHol{V}{U'}.
\eqno{\qed}
$$
\end{Corollary}

\begin{Example}\label{ex-sig-orientation}
For example, let $\sg = \sg (u,q)$ denote the expression
\begin{equation}\label{eq-sg-prod-expr}
    \sg = (\uhalf - \umhalf) 
              \prod_{n\geq 1} \frac{(1-q^{n}u) (1-q^{n}u^{-1})}
                                   {(1-q^{n})^{2}}.
\end{equation}
This may be considered as an element of
$\Z\psb{q}[u^{\pm\frac{1}{2}}]$ which is a  holomorphic function of
$(\uhalf,q)\in \C^{\times}\times D$, where $D = \{q\in \C |
0<\norm{q}<1 \}$. Let $\h = \{\tau\in \C | \Im \tau > 0 \}$ be the
open upper half 
plane.  We may consider $\sg$ as a holomorphic function of
$(z,\tau)\in \affan\times \h$ by
\begin{equation} \label{eq:42}
\begin{split}
      u^{r} & = e^{rz} \\
      q^{r}  & = e^{2\pi i r\tau}
\end{split}
\end{equation}
for $r\in \Q$.  It is easy to check using \eqref{eq-sg-prod-expr} that 
\begin{subequations} \label{eq-sigma-props}
\begin{align}
     \sigma (-z) & = -\sigma (z) \label{eq-sg-odd} \\
     \sigma (z)  & = z + o (z^{2}) \label{eq-sg-o-z}\\
     \sigma (uq^{n}) & = (-1)^{n} u^{-n} q^{-\frac{n^{2}}{2}} \sigma
(u). \label{eq-sg-theta}
\end{align}
\end{subequations}
The equations \eqref{eq-sigma-props} imply that the Taylor expansion
of $\sigma$ at the origin defines a
multiplicative analytic orientation
\begin{equation} \label{eq:20}
       MSO \xrightarrow{\Sigma} \HPCsym.
\end{equation}
\end{Example}

\begin{Definition} \label{def-sg-or}  
If $V$ is an oriented vector bundle, we write $\Sigma (V)$ for the
Thom class associated to the orientation \eqref{eq:20}, 
and  $\sigma (V)$ for the associated Euler class.
\end{Definition}

In \cite[\S 2.7]{AHS:ESWGTC} it is shown that $\Sigma$ is the sigma
orientation associated to the elliptic curve $C$: that is, the diagram 
\[
\xymatrix{
{MSO}
 \ar[r]^{\Sigma}
&
{\HPCsym}
\\
{\musix}
 \ar[u]
 \ar[ur]_{\sigma(\HPCsym,
C,\fmlgpof{\pr}^{-1})}
}
\]
commutes.

\section{Degree-four characteristic classes and theta functions}
\label{sec:degr-four-char}

\subsection{Degree-four characteristic classes} \label{sec:degr-four-char-1}

If $G$ is a connected compact Lie group, then by the splitting
principle the natural maps
\begin{equation} \label{eq:8}
\begin{split} 
    H^{2} (BG,\Z) & \rightarrow H^{2} (BT,\Z)^{W} \cong \chars^{W} \\
    H^{4} (BG,\Z) & \rightarrow  H^{4} (BT,\Z)^{W} \cong
                  (S^{2}\chars)^{W} \cong  \hom (\Gamma_{2}\cochars,\Z)^{W}
\end{split}
\end{equation}
are rational isomorphisms.  Here, if $M$ is an abelian group, then
$S^{2}M$ and $\Gamma_{2}M$ denote  
degree-two parts of the symmetric and divided power algebras on $M$.

Without rationalizing, a degree-four characteristic class $\dfcc \in
H^{4} (BG,\Z)$ gives rise to a homomorphism 
\[
    \Gamma_{2}\cochars \xra{I} \Z.
\]
We shall abuse notation and also write $I$ for the bilinear map 
\begin{equation} \label{eq:7}
    \cochars \times \cochars \xra{\gamma_{1}\times \gamma_{1}} 
    \Gamma_{2}\cochars \xrightarrow{I}\Z.
\end{equation}
We shall say that the characteristic class $\dfcc$ is 
\emph{positive 
definite} if the pairing $I$ is so.
We also write $\qfmsym$ for the quadratic function given by the composition 
\[
\xymatrix{
{\Gamma_{2}\cochars}
 \ar[r]^{I}
&
\Z
\\
{\cochars}
 \ar[u]^{\gamma_{2}}
 \ar[ur]_{\qfmsym}
}
\]
(If $A$ and $B$ are abelian groups, then a function 
\[
     f: A\to B
\]
is \emph{quadratic} if 
\begin{align*}
     0 & = f (0)\\
     0 & = f (x+y+z) - f (x+y) - f (x+z)- f (y+z) + f (x) + f (y) + f (z)\\
     f (-x) & = f (x).
\end{align*}
The function
\[
     A \xra{\gamma_{2}} \Gamma_{2} A
\]
is the universal quadratic function out of $A$).

From the definitions it follows that 
\begin{equation}\label{eq:13}
\begin{split}
    \qfm{a + b} & = \qfm{a} + I (a,b) + \qfm{b} \\
    \qfm{n a}   & = n^{2}\qfm{a} \\
    \qfm{w a }  & = \qfm{a } \\
    I (wa,wb) & = I (a,b)
\end{split}
\end{equation}
for $a, b\in \cochars$, $n\in \Z$,  and $w\in W$.

There are a variety of ways to express the relationship between the
characteristic class $\dfcc$ and the map $I$.  For example,
suppose that $Q_{0}$ and $Q_{1}$ are two principal $G$-bundles over $X$,
given as maps 
\[
    X \xrightarrow{Q_{i}} BT.
\]
Then we get a new principal $G$-bundle as the composition 
\[
        X \xrightarrow{\Delta} 
        X\times X  \xra{Q_{0}\times Q_{1}}
        BT\times BT 
        \xra{\mu}
        BT.
\]
The effect of $Q_{i}$ in cohomology 
\[
      \HZ{X} \xra{\HZ{(Q_{i})}} \HZ{BT}\cong \VG{\Ga}.
\]
is an $\HZ{X}$-valued point of $\VG{\Ga}$.  Equation~\eqref{eq:13}
implies that 
\[
     \dfcc (\mu (Q_{0}\times Q_{1})\Delta) = 
     \dfcc (Q_{0}) + I (\HZ{(Q_{0})}, \HZ{(Q_{1})}) + \dfcc (Q_{1}),
\]
where we have extended $I$ to a bilinear map 
\[
   (\VG{\Ga})\times (\VG{\Ga}) \rightarrow \Ga.
\]
As another example, 
suppose that $A \subseteq \T$ is a closed subgroup
with character group $\Z/n\Z$ ($n\geq 0$), and suppose that 
$m: A\to T$ is a  
homomorphism.  If 
\[
\barm, \barm' \in \hom (\T,T) = \cochars
\]
satisfy 
\[
    \barm\restr{A} = \barm'\restr{A} = m,
\]
then 
\[
    \barm' = \barm + n\delta
\]
for some $\delta\in \cochars$, and equation~\eqref{eq:13} implies that
\begin{align*}
   \qfm{\barm'} & \equiv \qfm{\barm} \mod n 
\end{align*}
We write $\qfm{m}$ for the class of $\qfm{\barm}$ in $\Z/n$.
If $z$ is the chosen generator of $H\Z^{2}B\T$, write also $z$
for the induced generator of $H\Z^{*}BA \cong (\Z/n\Z)[z]$.
With these conventions
\begin{equation} \label{eq:14}
    (Bm)^{*}\dfcc = \qfm{m} z^{2}.
\end{equation}
We write $\Ihatsym $ for the map 
\[
   \Ihatsym : \cochars \rightarrow \chars
\]
which is the adjoint of~\eqref{eq:7}.
Note that we have 
\[
\Ihat{w a} (w b) = I (wa,wb) = I (a,b) = \Ihat{a} (b).
\]
It follows that if $a \in \cochars^{W}$ then 
\[
   \Ihat{a} \in \chars^{W} \rightarrow H^{2} (BG,\Q).
\]
defines a rational characteristic class of principal $G$-bundles,
which we also denote $\Ihat{a}$.

Now suppose that $\PB$ is a $\T$-equivariant principal $G$-bundle over
a connected trivial $\T$-space $Y$, so in particular 
$\PB_{\T}$ is a 
principal $G$-bundle over $Y_{\T} = B\T\times Y$.   
Let $m \in \cochars$ be  a
reduction of the action of $\T$ on $\PB/Y$.   Then $\PB (m)$ is a
principal $Z (m)$-bundle over $Y$, and we have the following.

\begin{Lemma} \label{t-le-equiv-char-class-with-I-and-qfm}
In $H^{4} (Y_{\T},\Q) = H^{4} (B\T \times Y,\Q)$ we have
\[
\dfcc (\PB_{\T}) = \qfm{m}z^{2} + 
                  \Ihat{m} (\PB (m)) z + 
                  \dfcc (\PB).
\]
\end{Lemma}

\proof
By the splitting principle, we may suppose that we have a
factorization 
\[
\xymatrix{
&
{BT}
 \ar[d]
\\
{Y}
 \ar[ur]^{\PB}
 \ar[r]_{\PB}
&
{BG.}
}
\]
Lemma \ref{t-le-g-m} implies that the map 
\[
    B\T \times Y \xra{Bm \times Q} BT\times BT \xrightarrow{\mu} BT
\]
classifies $\PB_{\T}$.  It follows that 
\begin{align*}
    \dfcc (\PB_{\T}) & = \dfcc (\HQ{\mu (Bm\times\PB)})  \\
      & = 
    \dfcc (Bm) + I (\HQ{(Bm)}, \HQ{\PB}) + 
    \dfcc (\PB) \\
      & = \qfm{m} z^{2} + \Ihat{m} (\PB) z +  \dfcc (\PB).\tag*{\qed}
\end{align*}

\begin{Example}\label{ex-dfcc-spin}
Let $T_{SO (2d)}\cong \T^{d}$ be the standard maximal torus in $SO
(2d)$ (the image under the map $U (d) \rightarrow SO (2d)$ of the
torus of diagonal matrices), and let $T$ be its preimage in $Spin
(2d)$.  If $m= (m_{1},\dots ,m_{d})\in \Z^{d}\cong (T_{SO
(2d)})^{\vee}$, then there is a lift 
in the diagram  
\[
\xymatrix{
&
{T} 
 \ar[d]
\\
{\T}
 \ar[r]_-{m}
 \ar@{-->}[ur]
&
{T_{SO (2d)}}
}
\]
precisely when $\sum m_{i}$ is even, that is
\[
    \cochars \cong \{(m_{1},\dots ,m_{d})\in \Z^{d}| \sum m_{i} \text{
even} \}.
\]  
The function
\begin{equation} \label{eq:27}
     \qfm{m_{1},\dots ,m_{d}} = \frac{1}{2}\sum m_{i}^{2}
\end{equation}
therefore defines a quadratic map 
\[
    \qfmsym: \cochars\to \Z
\]
with associated bilinear form 
\[
    I (m,m') = \sum m_{i} m_{i}'.
\]
It is not hard to check using~\eqref{eq:27} that it is the quadratic
form associated to $c_{2}\in H\Z^{4} (BSpin (2d))$.

Now suppose that  
\[
   V = L_{1} \oplus \dots \oplus L_{d}
\]
is a $\T$-equivariant $Spin (2d)$ bundle,  written as a sum of
$\T$-equivariant complex line bundles, over a trivial $\T$-space $Y$.
Let $x_{i} = c_{1}L_{i}$, and suppose that $\T$ acts on $L_{i}$ by the
character $m_{i}$.  In order that $\borel{V}$ be a spin bundle, we
must have 
\[
    0 =  w_{2} (\borel{V}) =  \sum_{i} (m_{i} z + x_{i}) \pmod{2}.
\]
In that case, Lemma \ref{t-le-equiv-char-class-with-I-and-qfm} says
that 
\[
    c_{2} ( \borel{V}) = \left(\frac{1}{2} \sum m_{i}^{2}\right) z^{2} + 
                      (\sum m_{i} x_{i}) z + 
                      c_{2} V
\]
in $H^{4} (Y\times B\T)$. 
\end{Example}

\subsection{Theta functions}
\label{sec:theta-functions}

Recall from \S\ref{sec:elliptic-curves} that $q = e^{2\pi i \tau}$,
$\Lambda = 2\pi i \Z + 2 \pi i \tau \Z$, and that $\CQ$ is the elliptic curve  
\[
    \CQ = \Gaan/\Lambda \cong \Gman/q^{\Z}
\]
Following \cite{Looijenga:RootSystems}, we define a line bundle
over $\VG{\CQ}$ by the formula 
\begin{equation} \label{eq:21}
   \Loo =  \Loo (\dfcc) \eqdef \frac{(\VG{\Gman}) \times \C}
               {(u,\lambda) \sim 
                (uq^{m}, u^{\Ihat{m}}q^{\qfm{m}} \lambda)}.
\end{equation}
for $m\in \cochars$.

\begin{Remark}  \label{rem-2}
The identity map of $(\VG{\Gman})\times \C$ induces for $w\in W$ an
isomorphism of line bundles 
\[
    \Loo \cong w^{*}\Loo
\]
over $\VG{\CQ}$,
which is certainly compatible with the multiplication in
$W$, so $\Loo$ descends to a line bundle $\anomaly{\dfcc}$ on
$\VGW{\CQ}$.  
\end{Remark}

A theta function for $G$ of level $\dfcc$ is a
$W$-invariant holomorphic section of $\Loo (\dfcc)$, and so a section
of $\anomaly{\dfcc}$.

\begin{Definition}  \label{def-theta}
A \emph{theta function} for $G$ of level $\dfcc$
is a function 
\begin{equation} \label{eq:9}
   \theta  = \sum_{n>>-\infty} a_{n} (u) q^{n} \in (\Z[\chars])\lsb{q}
\end{equation}
which for $u=e^{z}$ and $q=e^{2\pi i \tau}$
is a holomorphic function of $(z,\tau)\in \VG{\affan}\times
\halfplane$, and which 
satisfies 
\begin{subequations}
\label{eq:10}
\begin{align}
     \theta (u q^{m} ) & = u^{-\Ihat{m}} q^{-\qfm{m}} \theta (u) \\
     \theta (u^{w}) & = \theta (u)
\end{align}
\end{subequations}
for $m\in \cochars$ and $w\in W$, where $\qfmsym$ and $I$ are the
quadratic form and bilinear map associated to the characteristic class
$\dfcc$.  
\end{Definition}

\begin{Remark} 
Lemma \ref{t-le-hol-char-classes-I} implies that a theta function for
$G$ gives a holomorphic characteristic class for principal $G$-bundles.
\end{Remark}

\begin{Remark} 
There is a good deal of redundancy in the definition.  Looijenga
studies \emph{formal} series of the
form~\eqref{eq:9} which transform according to \eqref{eq:10}.   
One has to be careful to identify a group formal series which is
closed under the operations implied by \eqref{eq:10}.  If $\dfcc$ is
positive definite, then every such formal theta function defines a
holomorphic function of $(z,\tau)$  (see \cite{Looijenga:RootSystems}).
\end{Remark}

\subsection{The sigma function and representations of
$\Loops{Spin (2d)}$} \label{sec:sigma-function-basic}

If $G$ is a simple and simply connected Lie group, then there
is a unique generator $\dfcc$ of $H^{4} (BG;\Z)\cong \Z$ such that the
associated pairing $I$ is positive definite.  If $\mathcal{V}$ is a 
representation of $\Loops{G}$  of level $k$ in the sense of
\cite{PressleySegal:Loo}, then its character $\chi$ is a theta
function for $G$ of level $k\dfcc$ \cite{Kac:Book}. 
The most important example for us is the representation
of $LSpin (2d)$ whose character is the Euler class of the sigma
orientation~\eqref{def-sg-or}; it is a theta function  
associated to the characteristic class $c_{2}$ of $Spin (2d)$. 

It is useful to be more explicit about this Euler class.  As in
Example~\ref{ex-dfcc-spin}, let $T_{SO (2d)}$ be the image
under the map $U (d) \rightarrow SO (2d)$ of the torus of diagonal
matrices.  For $u = (u_{1},\dots ,u_{d})\in T$ write
\begin{equation} \label{eq:18}
   \sigma_{d} (u) = \prod_{i=1}^{d}  \sigma (u_{i}).
\end{equation}
The product expression~\eqref{eq:18} and the fact \eqref{eq-sg-odd} that
$\sigma$ is odd imply that $\sigma_{d} (u^{w}) = \sigma_{d} (u)$ for
$w\in W$, and so $\sigma_{d}$ is a $W$-invariant function on
$\VG{\Gaan}$ (with zeroes precisely at the points $\VG{\Lambda}$).  
By Lemma~\ref{t-le-hol-char-classes-I} it defines a holomorphic
characteristic class for 
oriented vector bundles of rank $2d$.  If $V$ is such a vector bundle,
then $\sigma (V) = \sigma_{d} (V).$

The fractional powers of $u$   in the
expression~\eqref{eq-sg-prod-expr} for $\sigma$ prevent 
$\sigma_{d}$ from being a theta function for $SO (2d)$, but if $G=Spin
(2d)$ then the formula  
\[
    \sigma_{d} (u) = (-1)^{d} (\prod_{i} u_{i})^{\frac{1}{2}} 
                      \prod_{i} \left(
                        (1-u_{i})
                        \prod_{n\geq 1} 
                        \frac{( 1-q^{n} u_{i}) (1-q^{n}u_{i}^{-1})}
                             {(1-q^{n})^{2}} \right)
\]
shows that $\sigma_{d} \in (\Z[\chars]) \psb{q}$.  The
formula~\eqref{eq-sg-theta} implies that,
if $I$ and $\qfmsym$ are the pairing and quadratic
form associated to the generator $c_{2}\in H^{4} (BG;\Z)$ as in
Example~\ref{ex-dfcc-spin}, then  
\[
       \sigma_{d} (uq^{m}) = u^{-\Ihat{m}} q^{-\qfm{m}} \sigma_{d} (u),
\]
so $\sigma_{d}$ is a theta function for $Spin (2d)$ of level $c_{2}$.
Up to the factor $\prod_{n} (1-q^{n})^{2d}$, it is the
character of an irreducible representation of $\Loops{Spin
(2d)}$ of level $1$ \cite{Kac:Book,PressleySegal:Loo,Liu-1}.

\subsection{A useful holomorphic characteristic class}
\label{sec-useful-class}

Suppose that $\PB$ is a $\T$-equivariant principal $G$-bundle over a
connected $\T$-space $Y$.  Suppose that $\T[n]\subset \T$ acts
trivially on $Y$.  Let 
\[
  m: \T[n]\to T
\]
be a reduction of the action of $\T[n]$ on $\PB/Y$.

Let $\xi\in H^{4} (BG;\Z)$ be a positive definite class, with
associated quadratic form $\qfmsym $ and bilinear map $I$.   Let
$\theta$  be a theta function of level $\dfcc$ for $G$ (Definition
\ref{def-theta}).  Recall that 
\[
   \CQ = \Gaan / (2\pi i \Z + 2 \pi i \tau \Z), 
\]
and let $a$ be a point of $C$ of order $n$.
Choose  a point
$\bar{a}\in \Gaan$ such that $\pr (\bar{a}) = a$.  

We are going to define a
holomorphic function 
\[
F = F (\theta,m,\bara)  \in \O{} (\VG{\Gaan})^{W (m)},
\]
and so a holomorphic characteristic class of principal $Z (m)$-bundles
(see \S\ref{sec:holom-char-class}).
To give a formula for $F$ it is convenient to define
\begin{align*}     
\etorsa     & = e^{2 \pi i \bar{a}}, \\
\intertext{and recall \eqref{eq:43} that we have set}
     q^{r} & = e^{2\pi i r \tau} \\
     u^{r} & = e^{ r z} 
\end{align*}
for $r\in \Q$.

Since $na = 0$ in $C$, there are unique integers $\ell$ and $k$ such
that 
\[
      n\bar{a} = 2\pi i \ell + 2\pi i \tau k.
\]
We choose an extension $\bar{m}$ making the
diagram 
\[
\xymatrix{
{\T}
 \ar@{-->}[r]^{\bar{m}}
&
{T} \\
{\T[n]}
 \ar@{>->}[u]
 \ar[ur]_{m}
}
\]
commute.  With these choices, the formula for $F$ is 
\begin{align} \label{eq:23}
 F (z) & = u^{\frac{k}{n}\Ihat{\barm}}
         \etorsa^{\frac{k}{n}\qfm{\barm}}
         \theta (u \etorsa^{\barm}). 
\end{align}

\begin{Remark} \label{rem-3} 
The factors preceding the $\theta$ are closely related to the line
bundle $\mathcal{V}^{\frac{1}{n}}$ which appears in \cite{BottTaubes:Rig} and
\cite{AndoBasterra:WGEEC}.
\end{Remark}

\begin{Lemma}  \label{t-le-F-indep-m}
$F$ is independent of the choice of lift $\barm$.
\end{Lemma}

\proof
Let $\barm'$ be another choice.  Let $F'$ be the function defined
using $\barm'$.  Since $\barm'$ and $\barm$ both restrict to $m$
on $A$, there is a $\Delta\in \cochars$ such that 
\[
    \barm' = \barm + n \Delta.  
\]
We have 
\begin{align*}
  F' (z) & = u^{\frac{k}{n}\Ihat{\barm'}} 
             \etorsa^{\frac{k}{n} \qfm{\barm'}}
             \theta (u \etorsa^{\barm'}) \\
         & = u^{\frac{k}{n}\Ihat{\barm + n \Delta}}
         \etorsa^{\frac{k}{n} \qfm{\barm + n\Delta}}
              \theta (u \etorsa^{\barm} q^{k\Delta}) \\
         & = u^{\frac{k}{n} \Ihat{\barm}}
             u^{k\Ihat{\Delta}}
             \etorsa^{\frac{k}{n}\qfm{\barm}}
             \etorsa^{k I (\barm,\Delta)}
             q^{k^{2}\qfm{\Delta}}
             u^{-k\Ihat{\Delta}}
             \etorsa^{-k I (\Delta,\barm)}
             q^{-k^{2} \qfm{\Delta}}
             \theta (u\etorsa^{\barm}) \\
          & = F (z).\tag*{\qed}
\end{align*}

\begin{Lemma} \label{t-le-F-wm-F-m}
$F$ is invariant under $W (m)$.
\end{Lemma}

\begin{proof}
Suppose $w \in W (m)$.  We have 
\[
   w m = m
\]
so 
\[
   \barm = w \barm + n \Delta 
\]
for some $\Delta\in \cochars$.  The proof is now similar to the proof
of Lemma \ref{t-le-F-indep-m}.
\end{proof}

The dependence of $F (\theta,m,\bar{a})$ on the lift $\bara$ can be
calculated as follows.    Let $\bara'$ be another lift.  Then 
\[
    \bara'  = \bar{a} + \lambda
\]
for some $\lambda\in 2\pi i  \Z + 2 \pi i \tau \Z$, so letting
\[
   \etorsap = e^{\bara'}
\]
we have 
\[
   \etorsap = \etorsa q^{\delta}
\]
for some $\delta\in \Z$.  Thus 
\[
  \etorsap^{n} = q^{k'}
\]
with $k' = k + n \delta$, and so the quantity 
\[
    w (a,q^{\frac{1}{n}}) \eqdef \etorsa^{-1} q^{\frac{k}{n}} =
\etorsap^{-1}q^{\frac{k'}{n}}
\]
is an $n^{\text{th}}$ root of unity which is independent of the choice
of lift of $a$; in fact it is the Weil pairing of $a$ and $q^{\frac{1}{n}}$
in the curve $C$.  Because it is an $n^{\text{th}}$ root of unity, the quantity
\[
w (a,q^{\frac{1}{n}})^{\qfm{m}} \eqdef w (a,q^{\frac{1}{n}})^{\qfm{\barm}}
\]
is independent of the lift $\barm$ of $m$.  Using the functional
equation for $\theta$~\eqref{eq:10} it is straightforward to check the
following.

\begin{Lemma} \label{t-F-bara-weil}
$$
F (\theta,m,\bara') = w(a,q^{\frac{1}{n}})^{\delta\qfm{m}}
                        F (\theta,m,\bara).
\eqno{\qed}$$
\end{Lemma}

Lemma \ref{t-le-F-wm-F-m} implies that the Taylor series expansion of
$F (\theta,m,\bara)$ defines a class in  
$\HPsym (BZ (m)) = \HPsym (BT)^{W (m)}$, which we also denote $F
(\theta,m,\bara)$.   Let $\theta
(\PB,m,\bara) \in \HHol{Y}{\affan}$ be the holomorphic  
cohomology class given by the formula 
\begin{equation}\label{eq-theta-pb-m-bara}
     \theta (\PB,m,\bara) = \borel{\PB (m)}^{*} F (\theta,m,\bara)
\end{equation}
(using Lemma \ref{t-le-hol-char-classes-I} to conclude that the class
$\theta (\PB,m,\bara)$ is in fact holomorphic).

\begin{Lemma} \label{t-le-theta-pb-m-bara-indep-m}
If the centralizer $Z (m)$ is connected, for example if $G$ is unitary
or spin, then the class $\theta (\PB,m,\bara)$ is independent of the
reduction $m$ of the action of $\T[n]$ on $\PB/Y$. 
\end{Lemma}

\proof
Let $m'$ be another reduction of the action of $\T[n]$ on $\PB/Y$.  By
Lemma \ref{t-le-g-in-NT-if-Z-conn}, there is an element $w\in W$ such
that 
\[
   m' = w m.
\]
If $\barm: \T\to T$ is a lift of $m$, then $w \barm$ is a lift of
$m'$.  Using this lift to define $F (\theta,m',\bara)$, we have
\begin{align*}
     w^{*} F (\theta,m',\bara) (z)  & = F (\theta,m',\bara) (w (z)) \\
     & =  (u^{w})^{\frac{k}{n}\Ihat{w\barm}} 
          \etorsa^{\frac{k}{n}\qfm{w\barm}} 
          \theta (u^{w}\etorsa^{w\barm}) \\
     & = u^{\frac{k}{n}\Ihat{\barm}}
         \etorsa^{\frac{k}{n}\qfm{\barm}}
         \theta (u\etorsa^{\barm}) \\
     & = F (\theta,m,\bara) (z).
\end{align*}
Since the diagram 
\[
\xymatrix{
 {\borel{Y}}
 \ar[r]^-{\borel{Q (m)}}
 \ar[dr]_{\borel{Q (m')}}
& 
{BZ (m)}
 \ar[d]^{w}
\\
&
{BZ (m')}
}
\]
commutes, we have 
\begin{align*}
    \theta (\PB,m',\bara) & = \borel{Q (m')}^{*} F (\theta,m',\bara)  \\
                          & = \borel{Q (m)}^{*} w^{*} F (\theta,m',\bara) \\
                          & = \borel{Q (m)}^{*} F (\theta,m,\bara) \\
                          & = \theta (\PB,m,\bara).\tag*{\qed}
\end{align*}

The results of this section justify the following.  

\begin{Definition}\label{def-theta-PB-bara} 
Suppose that $Q$
is a $\T$-equivariant principal $G$-bundle over a connected $\T$-space
$Y$ on which $\T[n]$ acts trivially.  Suppose that for some
(equivalently every) reduction 
\[
    m: \T[n]\to T
\]
of the action of $\T[n]$ on $\PB/Y$, the centralizer $Z (m)$ is
connected.  Let $C$ be the elliptic curve
$\Gaan/2\pi i \Z + 2 \pi i \tau \Z$.  Let $a$ be a point of $C$ of order $n$.
Let $\theta$ be a theta function for $G$ of level
$\dfcc$, and let $\bara$ be a point of
$\affan$ whose image in $C$ is $a$.   We define $\theta (\PB,\bara) \in \HHol{Y}{\affan}$ to be the holomorphic
cohomology class 
\begin{equation} \label{eq:40}
    \theta (\PB,\bara) = \theta (\PB,m,\bara);
\end{equation}
in particular this class does not depend on the choice of $m$.
\end{Definition}

\begin{Lemma}  \label{t-le-theta-bara-dep}
If  $\bara' = \bara + 2\pi i s + 2 \pi i \delta \tau$
is another lift of $a$, then 
\begin{equation} \label{eq:41}
    \theta (\PB,\bara') = w (a,q^{\frac{1}{n}})^{\delta\qfm{m}} 
                          \theta (\PB,\bara);
\end{equation}
again this class does not depend on the choice of $m$. \qed 
\end{Lemma}

An important case of the preceding constructions is that $G=Spin
(2d)$, and $\sigma_{d}$ is the 
character of the basic representation of $\Loops{G}$ as in
\S\ref{sec:sigma-function-basic}.  Let $p : G\to SO (2d)$ be standard
the double cover.   Let $P/Y$ be 
the resulting $\T$-equivariant $SO (2d)$-bundle over $Y$, and let
$V$ be the associated vector bundle.  
We recall from \cite{BottTaubes:Rig} that Lemma \ref{t-le-Z-connected}
implies that  the sub-vector bundle $V^{\T[n]}$ is
orientable. Explicitly, if $m$ is a reduction of the action of $\T[n]$
on $\PB/Y$, then $pm$ is a reduction of the action of $\T[n]$ on
$P/Y$.   If $V^{\T[n]}$ has rank $2k$, then the map 
\[
    Y\xra{} BO (2k)
\]
classifying $V^{\T[n]}$ factors through the map $Y\to BZ (pm)$
classifying $P (pm)$, and we have the solid diagram 
\[
\xymatrix{
{BZ (m)}
 \ar[dd]_{Bp}
 \ar@{-->}[rr]
&&
{BSO (2k)}
 \ar[dd]
\\
&{Y} 
 \ar[ul]^{Q (m)}
 \ar[dr]^{V^{\T[n]}}
 \ar[dl]^{P (pm)}
\\
{BZ (pm)}
 \ar[rr]
&
&
{BO (2k)}
}
\]
Since (by Lemma \ref{t-le-Z-connected}) the centralizer $Z (m)$ is
connected, there is a dotted arrow making the diagram commute.  In other
words, $V^{\T[n]}$ is orientable, and a choice of orientation 
$\orient (V^{\T[n]})$ determines a map 
\[
    f: BZ (m) \rightarrow BSO (k)
\]
such that $fQ (m)$ classifies $(V^{\T[n]}, \orient (V^{\T[n]}))$.
In other words, Lemma \ref{t-le-Z-connected} implies the following
result (as I learned from \cite{BottTaubes:Rig}).

\begin{Lemma} \label{t-le-Z-connected-oriented}
If $V$ is a $\T$-equivariant $\spin (2d)$-vector bundle, then it is
$\T$-orientable.  \qed
\end{Lemma}

Let $\sigma_{k}\in \HPsym (BSO (2k))$ be the characteristic class
associated to the sigma function as in
\S\ref{sec:sigma-function-basic}.  We then have two $W (m)$-invariant
holomorphic functions on $\VG{\Gaan}$.  One is
$f^{*}\sigma_{k}$, and the other is $F (\theta,m,\bara)$.

\begin{Lemma} \label{t-le-R}
The ratio 
\[
    R = \frac{f^{*}\sigma_{k}}{ F (\sigma_{d},m,\bara)}
\]
is a $W (m)$-invariant unit of $\stalk{(\O{\VG{\Gaan}})}$.
\end{Lemma}

\begin{proof}
The poles of $R$ occur at zeroes of $F.$  Using the standard
maximal torus of $SO (2d)$ we may write a typical element of
$\VG{\Gman}$ as 
\[
     (u_{1},\dots ,u_{d}) \in (\Gman )^{d}\cong \VG{\Gman}.
\]
In these terms, a lift $\barm$ of $m$ is of the form
\[
     u^{\barm} = (u^{m_{1}},\dots ,u^{m_{d}})
\]
for some integers $m_{1},\dots ,m_{d}$.  
We have 
\[
       F (\sigma_{d},m,\bara) = 
       u^{\frac{k}{n}\Ihat{\barm}} \etorsa^{\frac{k}{n}\qfm{\barm}}
       \prod_{j=1}^{d} \sigma (u_{j} \etorsa^{m_{j}}).
\]
The product of sigma functions contributes a zero near $z=0$ if and
only if $m_{j}a = 0$ 
in $C$.  Let $j_{1},\dots ,j_{k}\in \{1,\dots ,d \}$ be the indices
such that $m_{j_{i}} a = 0$ in $C$; then 
\[
   f^{*}\sigma_{k} = \prod_{i=1}^{k}\sigma (u_{j_{i}} \etorsa^{m_{j_{i}}}).
\]
So the zeroes of $f^{*}\sigma_{k}$ precisely cancel those of 
$F (\sigma_{d}, m, \bara)$.
\end{proof}

Lemmas \ref{t-le-hol-char-classes-I} and \ref{t-le-R} together imply
that $R$  defines a holomorphic characteristic class for $Z
(m)$-bundles.  

\begin{Corollary} \label{t-co-ratio-euler-hol-nonvanish-finite}
The holomorphic characteristic class  
\[
    R (V,V^{\T[n]},\orient (V^{\T[n]}),\bara) \eqdef  \borel{Q
 (m)}^{*}R \in     (\HHstalk{Y})^{\times}
\]
is independent of the reduction $m$, and satisfies 
$$
    e_{\sigma} (\borel{V^{\T[n]}},\orient (V^{\T[n]})) = 
        R (V,V^{\T[n]},\orient (V^{\T[n]}),\bara)
        \sigma (V,\bara).
\eqno{\qed}$$
\end{Corollary}


\section{Equivariant elliptic cohomology}
\label{sec-equiv-ell}


\subsection{Adapted open cover of an elliptic curve}

If  $X$ is a $\T$-space and if $a$ is a point of $C$, then we define 
\[
    X^{a} = \begin{cases}
    X^{\circle[k]} & a \text{ is of order exactly }k\text{ in }C \\
    X^{\circle} & \text{otherwise}.
\end{cases}
\]
Let $N\geq 1$ be an integer.

\begin{Definition} \label{def-special}
A point $a \in C$ is \emph{special} for $X$ if $X^{a}\neq
X^{\circle}$.  
\end{Definition}

If $V$ is a $\circle$-bundle over a $\circle$-space $X$, then
it is convenient to consider a few additional points to be special.
Suppose that $F$ is a component of $X^{\T}$ and
\[
   m : \T \to T
\]
is a reduction of the action of $\T$ on the principal bundle
associated to $V$.  If we choose an isomorphism 
\[
     \cochars \cong \Z^{r},
\]
then we may view $m$ as an array of integers $(m_{1},m_{2},\dots
,m_{r})$.  These integers are called 
the \emph{exponents} or \emph{rotation numbers} of $V$ at $F$.  
Let $\compact{V}$ denote the one-point compactification of $V$.

\begin{Definition} \label{def-special-vector-bundle}
A point $a$ in $C$ is \emph{special} for $V$ if it is
special for $\compact{V}$ or if for some component $F$ of $X^{\T}$ there is a
rotation number $m_{j}$ of $V$ such that $m_{j}a = 0$.  
\end{Definition}

In either case, if $X$ is a finite $\T$-CW
complex, then the set of special points is a finite subset of the
torsion subgroup of $\CQ$.   

\begin{Definition} \label{def-adapted}
An indexed open cover $\{U_{a} \}_{a \in C}$ of $C$ is \emph{adapted to} $X$
or $V$ if
it satisfies the following.
\begin{enumerate}
\item [1)]  $a$ is contained in $U_{a}$ for all
$a\in C$.  
\item [2)] If $a$ is special and $a \neq b$, then $a\not\in U_{a}\cap U_{b}$.
\item [3)]  If $a$ and $b$ are both special and $a\neq b$, then the
intersection $U_{a}\cap U_{b}$ is empty.
\item [4)]  If $b$ is ordinary, then $U_{a}\cap U_{b}$ is non-empty
for at most one special $a$.
\item [5)]  Each $U_{a}$ is small~\eqref{def-small}.
\end{enumerate}
\end{Definition}

\begin{Lemma} \label{t-le-adapted-covers}
Let $X$ be a finite $\T$-CW complex.  Then 
$C$ has an adapted open cover, and any two adapted open covers
have a common refinement.  \qed
\end{Lemma}

\subsection{Complex elliptic cohomology}

Let 
\[
\EPsym \eqdef (\HPCsym,
C,\fmlgpof{\pr}^{-1}) = (\HPCsym,C,\log_{\omega})
\]
be the elliptic spectrum associated to the elliptic curve
$C$ (see \eqref{eq:19}).   Note that this is just a form of ordinary
cohomology.  The hat in the notation indicates that this is a
completion of the associated holomorphic cohomology.  Namely, 
suppose that $U\subset \CQ $ is a small open neighborhood of
the identity in $\CQ$.  Suppose that $V\subset \affan$ is the
component of $\pr^{-1}U$ containing the origin.  Let  
$\TE{X}{U} = (U,\Esym (\borel{X})\restr{U})$ be the ringed space defined as the pull-back in the diagram
\begin{equation} \label{eq:30}
\xymatrix{
{\TE{X}{U}} 
 \ar[r]
 \ar[d]
&
{\TH{X}\restr{V}}
 \ar[d] 
\\
{U}
 \ar[r]^{(\pr\restr{V})^{-1}}
&
{V.}
}
\end{equation}
The diagram~\eqref{eq:2} shows that $\TEP{X}$ and $\TE{X}{U}$
are related by the formula 
\[
    \TEP{X} \cong (\TE{X}{U})^{\wedge}_{0}.
\]

\subsection{Equivariant elliptic cohomology}

Grojnowski's circle-equivariant
extension of $\EPsym$ is a contravariant functor associating to 
a compact $\T$-manifold $X$ a $\Z/2$-graded $\O{C}$-algebra $\ETS{X}$, with the 
property that 
\[
    \ETS{\point} = \O{C}.
\]
It is equipped with a natural isomorphism 
\begin{equation} \label{eq:28}
    \EPsym (\borel{X})\xrightarrow[\cong]{\ac{X}} (\ETS{X})^{\wedge}_{0},
\end{equation}
such that 
\[
\ac{\point} = \log_{\omega} : (\O{C})^{\wedge}_{0} \cong \EPsym (B\T).
\]
We shall write $\ET{X}$ for the ringed space
$(C,\ETS{X})$ (see~\eqref{def-ringed-space}).  We take this
opportunity to phrase the account in 
\cite{AndoBasterra:WGEEC} of the construction of $\ETS{X}$  as the
construction of a covariant functor 
\[
   X \mapsto \ET{X}
\]
from finite $\T$-CW complexes to ringed spaces~\eqref{def-ringed-space}
over $C$, equipped with an identification 
\[
    \ET{\point} = C
\]
and a natural isomorphism of formal schemes 
\[
    \TEP{X} \xrightarrow[\cong]{\ac{X}} (\ET{X})^{\wedge}_{0}
\]
such that 
\[
\ac{\point}=\log_{\omega}:   \GpOf{\EPsym} = \Gah \cong \fmlgpof{C}.
\]
If $\mathcal{X} = (X,\O{X})$ is a
ringed space and $U$ is an open set of $X$, then we may write
$\mathcal{X} (U)$ in place of  $\O{X}(U)$.

Let $\{U_{a} \}_{a\in \CQ}$ be an adapted open cover of $\CQ$.  
For each $a\in \CQ$, we make a ringed space $\ETa{X}{a} = (U_{a},
\ETSa{X}{a})$ over 
$U_{a}$ as the pull back in the diagram
\begin{equation}\label{eq-et-local-sheaf}
\xymatrix{
{\ETa{X}{a}}
 \ar[r]
 \ar[d]
&
{\TE{X}{U_{a}-a}}
 \ar[d]
 \ar[r]
&
{\TH{X}\restr{V}}
 \ar[d]
\\
{U_{a}} 
 \ar[r]^{\trans_{-a}}
&
{U_{a}-a}
 \ar[r]^-{(\pr\restr{V})^{-1}} 
&
{V.}
}
\end{equation}
As in~\eqref{eq:30}, $V$ is the component of $\pr^{-1} (U_{a}-a)$
containing the origin.
In other words, let $V_{a}\subset \pr^{-1} (U_{a})$ be the component
containing the origin.  For $U\subset U_{a}$ let $V = V_{a} \cap
\pr^{-1} (U-a)$, and let
\[
    \ETSa{X}{a} (U) = \HHol{X^{a}}{V},
\]
considered as an $\O{C} (U)$-algebra via the isomorphism 
\[
   U \xrightarrow{\trans_{-a}} U-a \xra{(\pr\restr{V})^{-1}} V.
\]
If $a\neq b$ and $U_{a} \cap U_{b}$ is not empty, then by the definition
\eqref{def-adapted} of an adapted cover, at least one of $U_{a}$ and
$U_{b}$, suppose $U_{b}$, contains no special point.   In particular
we have $X^{b} = X^{\circle}$ and so an isomorphism
\begin{align*}\label{eq-fact-about-X-j}
   \TE{X^{b}}{U} & \cong \E{X^{b}} \times U \\
\intertext{i.e.}
   \Esym (\borel{X^{b}}) \tensor{\C[z]} \O{U} & \cong 
   \Esym (X^{b}) \tensor{\C} \O{U}
\end{align*}
for any small neighborhood $U$ of the origin.

\begin{Lemma} \label{t-le-restr-iso}
If $a\neq b$, $U\subset U_{a}\cap U_{b}$, and $b$ is not special,
then the  inclusion 
\[
    i: X^{b} \xra{} X^{a}
\]
induces an isomorphism
\begin{align*}
   \TE{X^{b}}{U-a}  \cong \TE{X^{a}}{U-a}.
\end{align*}
\end{Lemma}

\begin{proof}
If $a$ is not special, then $X^{a}=X^{b}$ and the  result is obvious.
If $a$ is special, then it is not contained in $U$ (by the definition
of an adapted cover), and so $0$ is not contained in $U-a$.    
The localization theorem \eqref{t-th-localization} gives the result.
\end{proof}

Let $U = U_{a} \cap U_{b}$.  We define 
\[
    \glue_{ab} = \glue_{ab}^{X}: \ETa{X}{a}\restr{U} \xra{\cong} 
                \ETa{X}{b}\restr{U}
\]
as the arrow making the diagram 
\begin{equation}\label{eq-phi-diagram}
\xymatrix{
{\ETa{X}{a}\restr{U}} 
 \ar[dr]
 \ar[rr]
 \ar[ddd]_{\glue_{ab}}
&
&
{\TE{X^{a}}{U - a}}
  \ar[d] 
& 
{\TE{X^{b}}{U - a}}
  \ar[dl]
  \ar[l]^{\cong}
  \ar[d]^{\cong}
\\
&
{U }
 \ar[r]^{\trans_{-a}} 
 \ar@{=}[d]
&
{U - a}
&
{\EP{X^{b}} \times (U - a)}
  \ar[l]
\\
&
{U}
 \ar[r]^{\trans_{-b}}
&
{U-b}
 \ar[u]_{\trans_{b-a}}
&
{\EP{X^{b}} \times (U - b)}
 \ar[l]
 \ar[u]_{\EP{X^{b}}\times \trans_{b-a}}
\\
{\ETa{X}{b}\restr{U}}
 \ar[rrr]
 \ar[ur]
&
&
&
{\TE{X^{b}}{U-b}}
 \ar[u]
 \ar[ul]
}
\end{equation}
commute.  The cocycle condition 
\[
    \glue_{bc} \glue_{ab} = \glue_{ac}
\]
needs to be checked only when two of $a,b,c$ are not special; and in that
case it follows easily from the equation 
\[
     \trans_{c-b} \trans_{b-a} = \trans_{c-a}.
\]
We shall write $\ET{X}$ for the ringed space over $C$, and $\ETS{X}$
for its structure sheaf.  One then has
the following (\cite{Grojnowski:Ell}; for a published account see \cite{Rosu:Rigidity}).

\begin{Proposition} \label{t-pr-ell-def}
$\ET{X}$ is a ringed space over $C$, which is independent up to canonical
isomorphism of the choice of adapted open cover. \qed
\end{Proposition}

\section{Equivariant elliptic cohomology of Thom spaces}
\label{sec:ellipt-cohom-thom-spaces}

Suppose that $V$ is a $\T$-equivariant vector bundle over $X$, and
let $V_{0}$ be the complement of the 
zero section of $V$.  We abbreviate as 
\begin{equation} \label{eq:33}
 \ETS{V}  = \ETS{V,V_{0}}
\end{equation}
the $\ETS{X}$-module associated to the reduced $\ETsym$-cohomology of
the Thom space of $V$.  Explicitly, for each point $a\in \CQ$ we
define a sheaf $\ETSa{V}{a}$ of $\ETSa{X}{a}$-modules as the
pull-back
\begin{equation}\label{eq-pull-back-for-V}
\xymatrix{
{  \ETSa{V}{a}}
    \ar[r]
    \ar[d] 
&
{\Esym (\borel{V^{a}})\restr{U_{a}-a}}
   \ar[d]
  \\
{   U_{a} }
    \ar[r]
  &
{   U_{a} - a.}
}
\end{equation}
For $a$ special and $b$ not with $U_{a}\cap U_{b}$ non-empty,
the isomorphism 
\[
   \glue_{ab}^{V}: \ETSa{V}{a}\restr{U_{a}\cap U_{b}} 
             \cong \ETSa{V}{b}\restr{U_{a}\cap U_{b}} 
\]
is given by 
\begin{multline*}
\Esym (\borel{V^{a}})\restr{U_{ab}-a} (U-a) \xrightarrow[\cong]{} 
\Esym (\borel{V^{b}})\restr{U_{ab}-a} (U-a) 
                     \cong \Esym (V^{b}) \tensor{\C} \O{} (U-a)
                     \xrightarrow{\trans_{b-a}^{*}}\\
                     \Esym (V^{b}) \tensor{\C} \O{} (U-b) \cong
                     \Esym (\borel{V^{b}})\restr{U_{ab}-b} (U-b);
\end{multline*}
where $U_{ab} = U_{a}\cap U_{b};$ $U$ denotes an open subset of
$U_{ab}$, and we have omitted a pullback \eqref{eq-pull-back-for-V} at either end.

\begin{Definition} \label{def-orientable}
The vector bundle $V$ is $\T$-\emph{orientable} if for each closed
subgroup $A\subseteq \T$, the fixed bundle $V^{A}$ over $X^{A}$ is
orientable.  A $\T$-\emph{orientation} $\orient$ on $V$ is a choice $\orient (V^{A})$
of orientation on $V^{A}$ for each $A$.
\end{Definition}

\begin{Lemma} \label{t-le-t-orientble-invertible}
If $V$ is a $\T$-orientable vector bundle over $X$, then $\ETS{V}$ is
a line bundle over $\ET{X}$, that is, $\ETS{V}$ is an invertible
$\ETS{X}$-module.
\end{Lemma}

\begin{proof}
If $V$ is $\T$-orientable, then the Thom isomorphism in ordinary
cohomology implies that each $\ETSa{V}{a}$ is an invertible
$\ETSa{X}{a}$-module.
\end{proof}

  We recall from 
\cite{Rosu:Rigidity,AndoBasterra:WGEEC} the construction of an 
explicit cocycle for this line bundle.
Let  $\phi$
be a multiplicative analytic orientation~\eqref{def-analytic-or}, and
let $\orient$ be a $\T$-orientation on $V$. 

\begin{Definition} \label{def-adapted-thom}
An
indexed open cover $\{U_{a} \}_{a\in C}$ of $C$ is 
\emph{adapted to} $(V,\phi,\orient)$ if it is
adapted to $V$ (see Definition~\ref{def-adapted}), and 
if for every point $a\in C$, the equivariant Thom class $\phi
(\borel{V^{a}},\orient)$ induces an isomorphism 
\[
\HHol{X}{U_{a}-a} \xrightarrow[\cong]{\phi}
\HHol{V^{a}}{U_{a}-a}.
\]
\end{Definition}

Corollary \ref{t-cor-analytic-thom} implies that $C$ has an indexed
open cover adapted to $(V,\phi,\orient)$.  Choose such a cover.  Suppose that
$a, b$ are two points of  $C$, such that $U = U_{a}\cap U_{b}$  is
non-empty: we may suppose 
that $a=b$ or that $b$ is not special.  Consider  the case that $b$ is
not special.  Let $\barU_{a}\subset \affan$ be the
component of the preimage of $U_{a}-a$ containing the origin, and let 
\[
   \pr\restr{\barU_{a}}: \barU_{a}\to U_{a}-a
\]
be the induced isomorphism.  Let $\barU \subset \barU_{a}$ be the preimage of
$U$.  Let  
\[
 j: (V^{b},V^{b}_{0}) \rightarrow (V^{a},V^{a}_{0})
\]
be the natural map.  The
Localization Theorem (Theorem \ref{t-th-localization}) implies that
the ratio of Euler classes
\[
   e_{\phi} (V^{a},V^{b},\orient) \eqdef \frac{j^{*}\phi (\borel{(
V^{a})},\orient (V^{a}))}{\phi (\borel{( V^{b})},\orient (V^{b}))} 
\]
is a unit of $\HHol{X^{\T}}{\barU}$.

Recall that there are tautological isomorphisms
\[
  \ETSa{X}{a} (U) 
       \xrightarrow[\cong]{\trans_{a}^{*}}
       \Esym (\borel{X^{a}})\restr{U-a} (U-a)
       \xrightarrow[\cong]{(\pr\restr{\barU})^{*}}
       \HHol{X^{a}}{\barU}.
\]
Let 
\[
e_{\phi} (a,b,\orient) \in \ETSa{X}{a} (U)^{\times}
\]
be given by the formula 
\begin{equation} \label{eq:6}
   (\pr\restr{\barU})^{*}\trans_{a}^{*}  e_{\phi} (a,b,\orient) = e_{\phi}
(V^{a},V^{b},\orient).
\end{equation}
Note that  $e_{\phi} (a,b,\orient) = 1$ if neither $a$ nor $b$ is special; we
also set 
\[
  e_{\phi} (a,a,\orient) = 1
\]
for all $a$.  It is easy to check that 
\begin{equation}\label{eq-cocycle-h-1}
      \glue_{bc}^{X} ( \glue_{ab}^{X} ( e_{\phi} (a,b,\orient))
e_{\phi} (b,c,\orient) )  = 
\glue_{ac}^{X} e_{\phi} (a,c,\orient)
\end{equation}
if $U_{a}\cap U_{b}\cap U_{c}$ is non-empty (since in that case at
least two of $a,b,c$ are ordinary).
Thus the $e_{\phi} (a,b)$ define a cohomology 
class $[\phi,V,\orient]\in H^{1} (C;\ETS{X}^{\times})$.  Let
$\ETS{X}^{[\phi,V,\orient]}$ be the resulting invertible sheaf of
$\ETS{ X}$-modules over $C$.  By construction we have the

\begin{Proposition} \label{t-pr-thom}
The Thom isomorphism $\phi$ induces an isomorphism 
\[
   \ETS{V} \cong \ETS{X}^{[\phi,V,\orient]}
\]
of $\ETS{X}$-modules. \qed
\end{Proposition}

In the
case of the orientation $\Sigma$ associated to the sigma
function~\eqref{def-sg-or}, we can be explicit about the open set on which 
\[
\sigma (V^{a},V^{b},\orient) \eqdef e_{\Sigma} (V^{a},V^{b},\orient)
\]
is a unit.

\begin{Lemma} \label{t-le-e-sigma-unit} 
Let $\Bar{B} \subset \affan$ be the preimage of the ordinary points of
$C$: that is, the complement of the (closed) set of points $\bara\in
\affan$ such that $X^{\pr (\bara)}\neq X^{\T}$ or $m \bara \in \Lambda$
for $m$ a character of the action of $\T$ on
$V\restr{X^{\T}}$.  Then 
\[
    \sigma (V^{\T[n]},V^{\T},\orient) \in \HHol{X^{\T}}{\Bar{B}}^{\times}.
\]
\end{Lemma}

\begin{proof}
Let $T$ be the standard maximal torus in $SO (2d)$, giving an
isomorphism 
\[
    \cochars \cong \Z^{d}.
\]
The reduction $m$ is then an array or integers $m= (m_{1},\dots
,m_{d})$.  It suffices to consider the case that 
\[
V\restr{X^{\T}}   \cong L_{1} + \dots  + L_{d}
\]
is a sum of line bundles, with $\T$ acting on $L_{i}$ by the character
$m_{i}$.  Let $x_{i} = c_{1}L_{i}$.   Then 
\begin{align*}
  \sigma (V^{\T[n]},V^{\T}) & = 
\frac{\prod_{m_{j}\equiv 0 \pmod{n}} \sigma (m_{j} z + x_{j})}
     {\prod_{m_{j} = 0} \sigma (m_{j}z + x_{j})}  \\
         & = \prod_{0\neq m_{j}\equiv 0 \pmod{n}} \sigma (m_{j}z + x_{j}).
\end{align*}
Since the $x_{j}$ are nilpotent, this is a unit in a neighborhood of
$z$ provided that $\prod_{0\neq m_{j}\equiv 0\pmod{n}}\sigma (m_{j}z )$ 
is non-zero.  This happens if and only if $m_{j}z\in \Lambda$.
\end{proof}

As mentioned in the introduction, if
f $W = (V_{0}, V_{1})$ is a pair of $\T$-oriented
vector bundles, then we extend the
notation~\eqref{eq:33} by defining 
$\ETS{W}$ to be the invertible sheaf
\[
   \ETS{W}  =  \ETS{V_{0}}\tensor{\ETS{X}}  \ETS{ V_{1}}^{-1}
\]
of $\ETS{X}$-modules.  As in Proposition \ref{t-pr-thom}, a choice of
$\T$-orientation $\Omega$ on $W$ gives rise to a class
$[\sigma,W,\Omega] \in H^{1}
(C,\ETS{X}^{\times})$,  
equipped with a canonical isomorphism of $\ETS{X}$-modules
\begin{equation} \label{eq:34}
    \ETS{W} \cong \ETS{X}^{[\sigma,W,\Omega]}.
\end{equation}

\section{A Thom class}\label{sec:thom-class}

\subsection{Equivariant elliptic cohomology of principal bundles}
\label{sec:equiv-ellipt-cohom-princ-I}

Suppose that $\PB/X$ is a $\T$-equivariant principal $G$-bundle over
$X$, so in particular $\borel{\PB}/\borel{X}$ is a principal
$G$-bundle as well.  Applying $\EPsym$ to the map 
\[
    \borel{X} \xra{} BG
\]
classifying $\borel{\PB}/\borel{X}$ induces a map 
\begin{equation} \label{eq:3}
     \TEP{X} \xrightarrow{\TEP{\PB}} \VG{\CQhat} / W
\end{equation}
as in \eqref{eq:47}.

If $F$ is a connected component of $X^{\T}$, and if 
\[
   m: \T\to T
\]
is a reduction of the action of $\T$ on $\PB\restr{F}$, 
then Lemma
\ref{t-le-equiv-char-class-restr-VC} implies that the diagram 
\begin{equation} \label{eq:4}
\xymatrix{
{\TEP{F}} 
 \ar@{=}[r]
 \ar[d]
&
{\EP{F}\times \CQhat}
 \ar[d]
 \ar[rr]^-{\TEP{\PB}}
& &
\VG{\CQ}/W
\\
{\ET{F}} 
 \ar@{=}[r]
&
{\EP{F}\times \CQ}
 \ar[rr]^-{\EP{\PB (m)} + m}
& & 
\VG{\CQ}/W (m)
 \ar[u]
}
\end{equation}
commutes.  It follows that if we define
\[
  \ET{\PB (m)} \eqdef \EP{\PB (m)} + m,
\]
then we have the commutative solid 
diagram 
\begin{equation}\label{eq-extend-equiv-cc-over-cq}
\xymatrix{
{\EP{F}\times \CQhat}
 \ar[rrr]^{\TEP{\PB (m)}}
 \ar[dr]
 \ar@{>->}[dd]
&
&
&
{\VG{\CQhat}/W (m)}
\ar@{>->}[ddd]
\ar[dl] \\
&
{\TEP{X}}
 \ar[r]^{\TEP{\PB}}
 \ar@{>->}[d] 
&
{\VG{\CQhat}/W}
 \ar@{>->}[d]
\\
{\EP{F}\times \CQ}
 \ar@{=}[d]
&
{\ET{X}}
 \ar@{-->}[r]_{\ET{\PB}}
&
{\VG{\CQ}/W}
\\
{\ET{F}}
 \ar[ur]
 \ar[rrr]_{\ET{\PB (m)}=\EP{\PB (m)}+m}
&
&
&
{\VG{\CQ}/W (m)}
 \ar[ul]
}
\end{equation}

In writing this paper, we were guided by the idea (Conjecture
\ref{t-co-princ-x-2}) that there should be a canonical map $\ET{\PB}$
making the whole diagram \eqref{eq-extend-equiv-cc-over-cq} commute.
We shall return to that 
question in \S\ref{sec:prin-bundles-ii}, discussing both why such a
map would be a good thing and why it is difficult to construct.  
For now we merely observe that the definition of $\ET{\PB (m)}$
implies the following. 

As we have already observed before Lemma
\ref{t-le-equiv-char-class-restr-VC},  the
addition  
\[
    (\VG{\CQ})\times (\VG{\CQ}) \rightarrow \VG{\CQ}
\]
induces a translation 
\[
    (\VG{\CQ})^{W (m)} \times \VG{\CQ}/W (m) \rightarrow \VG{\CQ}/W (m), 
\]
and so for $a\in C$ we get an operator
\[
   \trans_{a^{m}} : \VG{\CQ}/W (m) \rightarrow \VG{\CQ}/W (m).
\]
\begin{Lemma} \label{t-le-translation-on-F}
For $a\in \CQ$ the diagram 
\[
\begin{CD}
\EP{F}\times \CQ @> \ET{\PB (m)}    >> \VG{\CQ}/W (m) \\
@V \EP{F}\times \trans_{a} VV @VV \trans_{a^{m}} V \\
\EP{F}\times \CQ @> \ET{\PB (m)}    >> \VG{\CQ}/W (m)
\end{CD}
\]
commutes. \qed
\end{Lemma}

\subsection{The Thom class}

Let $G$ be a spinor group, and let $G'$ be a simple and simply
connected compact Lie group, a unitary group, or indeed any compact
connected Lie group with the property that the centralizer of any
element is connected.  Let $V$ be a 
$\T$-equivariant $G$-vector bundle over a finite $\T$-CW complex $X$,
and let $V'$ be a $\T$-equivariant $G'$-bundle (by which we mean the
vector bundle associated to a principal $G'$ bundle via a linear
representation of $G'$).  
Suppose 
that $\dfcc'$ is a degree-four characteristic classes for 
$G'$, with the property 
that 
\begin{equation} \label{eq:36}
   c_{2} (\borel{V}) = \dfcc' (\borel{(V')}).
\end{equation}
Suppose finally that $\theta'$ is a theta function for $G'$ of level
$\dfcc'$.  In this section we prove the following

\begin{Theorem} \label{t-th-thom-class}
A $\T$-orientation $\epsilon$ on $V$ determines a canonical global section
$\gamma=\gamma (V,V',\epsilon)$ of $\ETS{V}^{-1}$, such that
\begin{equation} \label{eq:32}
\gamma_{0} = \theta' (\borel{V'})\Sigma
(\borel{V})^{-1}
\end{equation}
under the isomorphism~\eqref{eq:28}
\[
(\ETS{V}^{-1})^{\wedge}_{0} \cong 
\EPsym (\borel{V})^{-1} = \HPsym (\borel{V})^{-1}
\]
of $\HPsym (\borel{X})$-modules. The formation of $\gamma$ is natural
in the sense if  $f: Z\to X$ is a map of finite $\T$-CW complexes, and 
\begin{align*}
     g/f: W/Z & \rightarrow V/X \\
     g'/f: W'/Z & \rightarrow V'/X
\end{align*}
are pull-backs of vector bundles, then under the isomorphism 
\begin{align}
    \ETS{g/f}: \ETS{f}^{*} \ETS{V}  & \cong \ETS{W}, \notag \\
\intertext{we have}
    \ETS{f}^{*} \gamma (V,V',\epsilon) 
    & = \gamma (W,W',g^{*}\epsilon). \label{eq:31}
\end{align}
\end{Theorem}

Let $\{U_{a} \}_{a\in C}$ be an indexed open cover  of $C$ adapted to
$(V,\Sigma,\orient)$ (Definition \ref{def-adapted-thom}).   By
Proposition \ref{t-pr-thom}, to give a section of $\ETS{V}^{-1}$ 
satisfying \eqref{eq:32}  is equivalent to giving a global
section of the sheaf 
$\ETS{X}^{-[\sigma,V,\orient]}$ whose value in 
\[
\stalk{(\ETS{X}^{-[\sigma,V,\orient]})}^{\wedge}\cong \HPsym (\borel{X})
\]
is $\theta' (\borel{V'})$; this is what we shall do.  The description
of $\ETS{X}^{-[\sigma,V,\orient]}$ in
\S\ref{sec:ellipt-cohom-thom-spaces} shows that this amounts to giving
sections $\gamma_{a} \in \ETSa{X}{a} (U_{a})$ for $a\in C$, such that
$\gamma_{0} = \theta' (\borel{V'})$ and 
\[
            \glue_{ab} (\sigma (a,b,\orient ) \gamma_{a}) = \gamma_{b}
\]
when $a$ is special and $b$ is ordinary.

Let $B\subset C$ be the set of ordinary points, and, as in Lemma
\ref{t-le-e-sigma-unit}, let $\Bar{B}\subset \affan$ be the preimage
of $B$ in $\affan$.  Lemma \ref{t-le-e-sigma-unit} tells us that the 
formula 
\[
   \Bar{\gamma}_{B} \eqdef \theta' (\borel{V'})\restr{X^{\T}}
                           \frac{\sigma (V^{\T})}
 {\sigma (\borel{V}\restr{X^{\T}})}
\]
defines an element of $\HHol{X^{\T}}{\Bar{B}}$.

\begin{Lemma} \label{t-le-gamma-B}
 For $\lambda\in \Lambda$, 
\[
       \trans_{\lambda}^{*} \Bar{\gamma}_{B} = \Bar{\gamma}_{B}.
\] 
\end{Lemma}

\begin{proof}
There is a $k$ such that 
\[
   q^{k} = e^{\lambda}.
\]
Let $\PB$ be the principal $G$-bundle associated to $V$, and let
$\PB'$ be the principal bundle associated to $G'$.  If $F$ is a
component of $X^{\T}$, let $m: \T\to T$ be a reduction of the action
of $\T$ on $\PB\restr{F}$, and let $m': \T\to T'$ be a reduction of 
the action of $\T$ on $\PB'\restr{F}$.   

The principal bundle $\PB (m)/F$ is classified by a map 
\[
    F \xra{\PB (m)} BZ (m);
\]
in rational cohomology this becomes 
\[
   \HQ{F} \xra{\HQ{\PB (m)}} \HQ{BZ (m)} \cong (\VG{(\Gah)_{\Q}})/W (m),
\]
an $\HQ{F}$-valued point of $(\VG{(\Gah)_{\Q}})/W (m)$.  Since $F$ has
the homotopy type of a finite CW-complex, the reduced cohomology of
$F$ is nilpotent, and so we may consider $\exp (\HQ{\PB (m)})$ as an
$\HQ{F}$-valued point of $(\VG{(\Gm)_{\Q}})/W (m).$  

Let
\begin{align*}
     D & = \exp (\HQ{\PB (m)})\in (\VG{\Gm}/W (m))(\HQ{F}) \\
     D' & = \exp (\HQ{\PB' (m')}) \in (\VG{\Gm}'/W' (m'))(\HQ{F}) \\
     u  & = \exp (z).
\end{align*}
Then 
\begin{equation} \label{eq:29}
       \sigma (\borel{V}\restr{F}) (z)  = \sigma_{d} (D u^{m})
\end{equation}
and 
\begin{align*}
      \trans_{\lambda}^{*}\sigma (\borel{V}) (z)&  = 
      \sigma_{d} (D u^{m}q^{km})  \\
      &  = 
     (D u^{m})^{-k\Ihat{m}} q^{-k^{2}\qfm{m}} \sigma (\borel{V}) (z) \\
      & = 
D^{-k\Ihat{m}} u^{-k I (m,m)} q^{-k^{2}\qfm{m}} \sigma (\borel{V}) (z).
\end{align*}
Similarly 
\[
     \trans_{\lambda}^{*}\theta' (\borel{V'}) (z)  = 
     (D')^{-k\Ihatsym' (m')} 
     u^{-k I' (m',m')} 
     q^{-k^{2} \qfmsym' (m')} 
     \theta' (\borel{V'}) (z).
\]
If $\borel{c_{2} (V)} = \borel{\dfcc' (V')}$, then Lemma
\ref{t-le-equiv-char-class-with-I-and-qfm} 
implies that  
\begin{align*}
            D^{\Ihat{m}} & = 
            \exp (\Ihat{m} (\HQ{\PB (m)})) = 
            \exp (\Ihat{m'} (\HQ{\PB' (m')})) = (D')^{\Ihatsym' (m')} \\
            \qfm{m} & = \qfmsym' (m')
\end{align*}
which gives the result.
\end{proof}

\begin{Example}
To illustrate the notation used in the proof, suppose we have chosen
an isomorphism $T\cong \T^{d}$ and so also $\cochars\cong \Z^{d}$.
Then we have 
\[
   \VG{\Gah} \cong \mathbb{\widehat{A}}^{d} \cong \spf
   \Z\psb{x_{1},\dots ,x_{d}}. 
\]
Suppose that the map $F\rightarrow BG$ classifying $V\restr{F}$
factors through $BT$, i.e. that
\[
   V\restr{F} = L_{1} \oplus \dots \oplus L_{d}
\]
is written as a sum of complex line bundles.  Then under the resulting
map 
\[
    \HQ{F} \rightarrow (\VG{(\Gah)_{\Q}}),
\]
the coordinate function $x_{j}$ pulls back to $c_{1}L_{j}$.  If $\T$
acts on $L_{j}$ by the character $m_{j}$, then equation~\eqref{eq:29}
becomes the more familiar equation
\[
\sigma (\borel{V}\restr{F}) (z)  = \prod_{j} \sigma (e^{x_{j} + m_{j}z}).
\] 
\end{Example}

Lemma \ref{t-le-gamma-B} says that $\Bar{\gamma}_{B}$ descends to a
function $\gamma_{B}$ on 
\[
    \E{X^{\T}}\times B \subset \E{X^{\T}}\times C.
\]
For $b$ an ordinary point of $C$, we define 
\[
    \gamma_{b} \eqdef \gamma_{B}\restr{U_{b}} \in \O{}
(\E{X^{b}}\times U_{b}) = \ETSa{X}{b} (U_{b}).
\]

Now suppose that $a$ is a special point.  Let $\bara$ be a preimage  of $a$
in $\affan$, and define $\gamma_{a} \in \ETSa{X}{a} (U_{a})$ by the
formula  
\[
   \gamma_{a} \eqdef  \trans_{-a}^{*}((\pr\restr{W})^{-1})^{*}
    ( R(V,V^{a},\orient (V^{a}),\bara) \theta' (V',\bara)),
\]
where $W\subset \pr^{-1}U_{a}$ is the component containing the origin
(see~\eqref{eq-et-local-sheaf}) and $R$ is the characteristic class
defined in Corollary \ref{t-co-ratio-euler-hol-nonvanish-finite}.
This is a definition in view of the  

\begin{Lemma} \label{t-le-gamma-a}
The class $\gamma_{a}$ is independent of the lift $\bara$.
\end{Lemma}

\begin{proof}
Suppose that $\bara$ and $\bara'$ are two lifts of $a$.  Let $\etorsa$, $B$,
and $\delta$ be given by 
\begin{align*}
    \etorsa & = \exp (\bara)  \\
    \etorsap & = \exp (\bara') \\
    \etorsap & = \etorsa q^{\delta}.
\end{align*}
Let $Y$ be a component of $X^{a}$. Let $m$ be a reduction of the
action on $\T$ on $\PB\restr{Y}$, and let $m'$ be a reduction of the
action of $\T$ on $\PB'\restr{Y}$.
Lemma~\ref{t-le-theta-bara-dep} implies that
\begin{align*}
      \theta' (\PB',\bara') & = w (a,q^{\frac{1}{n}})^{\delta\qfm{m}}
      \theta' (\PB',\bara) \\
      R (V,V^{a},\orient (V^{a}),\bara') & = 
                            w (a,q^{\frac{1}{n}})^{-\delta\qfm{m'}}
      R (V,V^{a},\orient (V^{a}),\bara).
\end{align*}
Equation~\eqref{eq:14} 
implies that 
\[
     \qfm{m'} \equiv \qfm{m} \mod n
\]
which gives the result.
\end{proof}

\begin{Lemma}  \label{t-le-gamma-glues}
The various sections $\gamma_{a}$ for $a\in C$ define a global section
of $\ETS{X}^{-[\sigma,V]}$, whose value in
$\stalk{(\ETS{X}^{-[\sigma,V]})}^{\wedge}$ is $\theta' (V')$.
\end{Lemma}

\begin{proof}
The value at $0$ follows from the fact that 
\[
   R (V,V,0) = 1, 
\]
as is easily checked.  To see that the $\gamma_{a}$ assemble into a
global section of $\ETS{X}^{-[\sigma,V]}$, we must show that, if $a$ is
a special point of order $n$ and $b$ is ordinary with $U = U_{a}\cap
U_{b}$ nonempty, then
\begin{equation} \label{eq:5}
    \glue_{ab} (\sigma (a,b,\orient ) \gamma_{ab}) = \gamma_{b}.
\end{equation}
Let $i: X^{b}\to X^{a}$ be the inclusion.  The diagram
\eqref{eq-phi-diagram} defining $\glue_{ab}$ together with the
equation \eqref{eq:6} for $\sigma (a,b,\orient)$ reduces 
\eqref{eq:5} 
to
\[
\trans_{b-a}^{*} 
\left(\sigma (V^{a},V^{b},\orient)^{-1} i^{*}\trans_{a}^{*}\gamma_{a}
\right) =  
\trans_{b}^{*}\gamma_{b}
\]
or equivalently 
\begin{equation} \label{eq:37}
\sigma (V^{a},V^{b},\orient)^{-1} i^{*}\trans_{a}^{*}\gamma_{a} =
\trans_{a}^{*} \gamma_{b}. 
\end{equation}
Suppose that $Y$ is a component of $X^{a}$, and suppose that $F$ is a
component of $X^{b}$ contained in $Y$.  Let $\barm$ be a
reduction of the action of $\T$ on $\PB\restr{F}$, and let $\barm'$ be a
reduction of the action of $\T$ on $\PB'\restr{F}$.     Lemma
\ref{t-le-equiv-char-class-with-I-and-qfm} and the 
equation $c_{2} (\borel{V}) = \dfcc' (\borel{(V')})$ imply that 
\begin{align}
     \qfm{\barm}  & = \qfmsym' (\barm') \label{eq:51} \\
     \Ihat{\barm} (\PB (\barm)) & = 
     \Ihatsym' (\barm') (\PB' (\barm')).\notag
\end{align}
As in the proof of Lemma \ref{t-le-gamma-B}, we let 
\begin{align*}
     D & = \exp (\HQ{\PB (\barm)})\in (\VG{\Gm}/W (m))(\HQ{F}) \\
     D' & = \exp (\HQ{\PB' (\barm')}) \in (\VG{\Gm}'/W' (m'))(\HQ{F}) \\
     u  & = \exp (z),
\end{align*}
and observe that Lemma \ref{t-le-equiv-char-class-with-I-and-qfm}
implies that 
\begin{equation} \label{eq:50}
   (Du^{\barm})^{\Ihat{\barm}} = (D'u^{\barm'})^{\Ihatsym (\barm')}.
\end{equation}
Let $\bara$ be a preimage of $a$; let $k$ be the integer such that 
\[
      n \bara = 2 \pi i \ell + 2\pi i \tau k;
\]
and let $\etorsa= \exp (\bara)$.

By Lemma \ref{t-le-m-Y-restr-to-Fcheck}, $\barm\restr{\T[n]}$ is a
reduction of the action of $\T[n]$ on $\PB\restr{Y}$ and similarly for
$\barm'$. By  
Corollary~\ref{t-co-ratio-euler-hol-nonvanish-finite} and
Lemma~\ref{t-le-F-indep-m}, 
we may use $\barm$ to calculate $R (V,V^{a},\orient (V^{a}),\bara)$
and $\barm'$ to calculate $\theta (\borel{V'},\bara)$.  The 
restriction to $F$ of the left side of \eqref{eq:37} becomes 
\begin{align*}
    \frac{\sigma (V^{b},\orient )}
         {\sigma (\borel{V^{a}},\orient )\restr{F}}
    R (V,V^{a},\orient,\bara)\restr{F} \theta' (\borel{V'},\bara) \restr{F} 
 = & \sigma (V^{b},\orient)  \\
   & \frac{(D'u^{\barm'})^{\frac{k}{n}\Ihatsym' (\barm')} 
         \etorsa^{\frac{k}{n}\qfmsym' (\barm')}
         (\trans_{\bara^{\barm'}}\theta') (\borel{V'}\restr{F})}
     {(Du^{\barm})^{\frac{k}{n}\Ihat{\barm}} 
         \etorsa^{\frac{k}{n}\qfm{\barm}}
          (\trans_{\bara^{\barm}}\sigma) (\borel{V}\restr{F})}
     \\
 = &\frac{\sigma (V^{b},\orient)}
         {(\trans_{\bara^{\barm}}\sigma) (\borel{V}\restr{F})}
        (\trans_{\bara^{\barm'}}\theta') (\borel{V'}\restr{F}) \\
= &
\trans_{a}^{*}
\left(\frac{\sigma (V^{b},\orient)}
           {\sigma (\borel{V}\restr{F})}  \theta' (\borel{V'}\restr{F})
\right) \\
= &  \trans_{a}^{*}\gamma_{b.}
\end{align*}
In the second equation we have equations~\eqref{eq:51}
and~\eqref{eq:50}.  In the third equation we
have used 
Lemma~\ref{t-le-translation-on-F} and
the fact that $\sigma (\borel{(V^{b})})= \sigma (V^{\T})$ is invariant under
translation.  
\end{proof}

This completes the construction of the section $\gamma$ promised in
Theorem \ref{t-th-thom-class}.  The naturality \eqref{eq:31} is
straightforward, given the canonical nature of the sections $\gamma_{a}$.

\section{The sigma orientation}
\label{sec:sigma-orientation}

Now suppose that $V= (V_{0},V_{1})$ is a pair of $\T$-equivariant spin
vector bundles of even rank over a finite $\T$-CW complex $X$, with
the property that  
\[
     c_{2} (\borel{V})  = 0.
\]
Recall (\cite{BS:atmss,BottTaubes:Rig}; see
Lemma \ref{t-le-Z-connected-oriented}) that each $V_{i}$ is
$\T$-orientable.   

\begin{Theorem} \label{t-th-sigma-orientation}
A $\T$-orientation $\orient$ on $V_{0}$ and $V_{1}$ determines a
canonical trivialization 
$\gamma (V) = \gamma (V,\orient)$ of $\ETS{V}$, whose 
value in $\ETS{V}^{\wedge}_{0} \cong \HPsym ((\borel{X})^{\borel{V}})$
is $\Sigma (\borel{V})$ (see Definition \ref{def-sg-or}). 
The association $V\mapsto \gamma (V)$ is stable, natural, and
exponential in the sense described in Theorem \ref{t-th-sigma-intro}.
\end{Theorem}

\begin{proof} The proof proceeds much as the proof of
Theorem~\ref{t-th-thom-class}: 
to construct $\gamma (V)$ it is equivalent to give a section $\gamma$
of $\ETS{X}^{[\sigma, V,\orient ]}$ whose value in
\[
(\ETS{X}^{[\sigma,V,\orient]})^{\wedge}_{0}\cong \HPsym (\borel{X})
\]
is $1$.  

Once again, let $B\subset C$ be the set of ordinary points, and, as in Lemma
\ref{t-le-e-sigma-unit}, let $\Bar{B}\subset \affan$ be the preimage
of $B$ in $\affan$.  Lemma \ref{t-le-e-sigma-unit} tells us that the 
formula 
\[
   \Bar{\gamma}_{B} \eqdef 
    \frac{\sigma (\borel{(V_{0})})}{\sigma (V_{0}^{\T})}
    \frac{\sigma (V_{1}^{\T})}{\sigma (\borel{(V_{1})})}
\]
defines an \emph{unit} in  $\HHol{X^{\T}}{\Bar{B}}^{\times}$.  The
same argument as 
in Lemma~\ref{t-le-gamma-B} shows once again that $\Bar{\gamma}_{B}$
descends to a function $\gamma_{B}$ on  
\[
    \E{X^{\T}}\times B \subset \E{X^{\T}}\times C.
\]
For $b$ an ordinary point of $C$, we define 
\[
    \gamma_{b} \eqdef \gamma_{B}\restr{U_{b}} \in \O{}
(\E{X^{b}}\times U_{b})^{\times} = \ETSa{X}{b} (U_{b})^{\times}. 
\]
Now suppose that $a$ is a special point.  Let $\bara$ be a preimage  of $a$
in $\affan$, and define $\gamma_{a}$ by the formula 
\[
    (\pr\restr{W})^{*}\trans_{a}^{*} \gamma_{a} \eqdef 
 \frac{R (V_{1},V_{1}^{a},\orient (V_{1}^{a}),\bara)}
      {R (V_{0},V_{0}^{a},\orient (V_{0}^{a}),\bara)},
\]
where $W\subset \pr^{-1} ( U_{a})$ is the component containing the origin, 
and $R$ is the characteristic class defined in Corollary
\ref{t-co-ratio-euler-hol-nonvanish-finite}. As in
Lemma~\ref{t-le-gamma-a}, $\gamma_{a}$ is independent of the 
lift $\bara$.  In this case, additionally, 
Corollary~\ref{t-co-ratio-euler-hol-nonvanish-finite} implies that
$\gamma_{a}$ is a \emph{unit}, i.e. an element of $\ETSa{X}{a}
(U_{a})^{\times}$.

The same argument as in the proof of Lemma \ref{t-le-gamma-glues}
shows that the sections $\gamma_{a}$ for $a\in C$ assemble into a global
section of $\ETS{X}^{[\sigma,V,\orient]}$, which is a
trivialization because it is so on each $U_{a}$.

The fact that the section $\gamma (V)$ is independent of the choice of
$V_{i}$, as well as the 
fact that $\gamma (V\oplus V') = \gamma (V) \otimes \gamma (V')$
follows from definition of $\gamma (V)$ and the equation 
\[
    \sigma (V\oplus V') = \sigma (V) \sigma (V').
\]
The naturality under change of base is clear from the construction.
\end{proof}

\section{A conceptual construction of the equivariant\nl sigma orientation}
\label{sec:prin-bundles-ii}

This section is devoted to a discussion of the following elaboration
of Conjecture \ref{t-co-princ-x}.

\begin{Conjecture}  \label{t-co-princ-x-2}
Equivariant elliptic cohomology ought to have the
following feature.  Suppose that $\PB$ is a $\T$-equivariant principal 
$G$-bundle over a $\T$-space 
$X$.  Then there is a canonical map 
\[
    \ET{X} \xra{\ET{\PB}} \VGW{\CQ}
\]
making the diagram 
\[
\begin{CD}
    \TEP{X} @> \TEP{\PB} >> \VGW{\fmlgpof{\CQ}} \\
    @VVV @VVV \\ 
    \ET{X} @> \ET{\PB} >> \VGW{\CQ}
\end{CD} 
\]
commute and having all the properties listed in Conjecture
\ref{t-co-princ-x}.   Moreover, 
for all components $F$ of $X^{\T}$ and all reductions 
\[
m: \T\to T
\]
of the action of $\T$ on $\PB\restr{F}$, 
the diagram \eqref{eq-extend-equiv-cc-over-cq} should commute.
\end{Conjecture}

\subsection{Why the conjecture should be true, and why it is
nevertheless difficult to prove}

As explained in the introduction, the conjecture gives an elegant
description of the equivariant sigma orientation, which even
illuminates the non-equivariant case. 

Although we have stated a conjecture, it is really a proposal 
of structure which should be sought in \emph{some}
equivariant elliptic cohomology theory.  How difficult it is to
establish the conjecture depends on your ontology.  Given
a fully developed theory of equivariant elliptic cohomology as
proposed by Ginzburg-Kapranov-Vasserot, it is not difficult at least
to construct the map $\ET{\PB}$.  

To see this, recall that Ginzburg, 
Kapranov, and Vasserot have proposed that equivariant elliptic
cohomology for the curve $C$ and the (compact connected Lie) group $G$
should be a covariant functor 
\[
     \Egkv{G}{\slot}  : \CategoryOf{$G$--spaces} \rightarrow 
                    \CategoryOf{schemes over $\VGW{C}$}.
\]
If $\PB$ is a $\T$-equivariant principal $G$-bundle over
$X$, then the $(\T\times G)$-equivariant elliptic cohomology of $\PB$
should be a scheme
\[
   \Egkv{\T\times G}{\PB}   \rightarrow C \times (\VG{C}/W)
\]
over both $C$ (via the $\T$-action) and $\VG{C}/W$ (via the
$G$-action).  Now $G$ acts freely on $\PB$ with quotient $X$, so one
expects that  
\[
   \Egkv{\T \times G}{\PB}  \cong \Egkv{\T}{X}.
\]
Combining these observations leads to the prediction that a
$\T$-equivariant principal $G$-bundle over $X$ should give rise to the
map 
\[
     \Egkv{\T}{X} \xrightarrow{\Egkv{\T}{\PB}}  \VGW{C}.
\]
There are two problems with this proposal.  First, we have in this
paper described only $\T$-equivariant elliptic cohomology.   It may
not be difficult to resolve this problem: it is not
so difficult to imagine an analogous construction of $G$-equivariant
elliptic cohomology, and indeed Grojnowski does so in
\cite{Grojnowski:Ell}.  That leaves the second problem, that 
the underlying space of Grojnowski's functor 
is always $C$: 
if 
\[
   \pi: \Egkv{\T}{X} \rightarrow C
\]
is the structural map associated to an elliptic cohomology which is proposed by
Ginzburg-Kapranov-Vasserot, then we have worked in this paper with the
sheaf 
\[
   \ETS{X}  = \pi_{*}\O{\Egkv{\T}{X}}.
\]
To see why this is a problem, note that in order to construct a map 
\[
    \ET{\PB}: \ET{X} \rightarrow \VG{C}/W
\]
we must in particular construct a map of topological spaces 
\begin{equation} \label{eq:15}
       C \rightarrow \VGW{C}.
\end{equation}
If $X^{\T}$ is non-empty, then $\ET{X^{T}} = X^{\T}\times C$, and
for each connected component $F$ of $X^{\T}$,  a reduction of the
action of $\T$ on  $\PB\restr{F}$ gives a map 
\[
    m: C\to \VG{C}
\]
and so one has a place to start.  But it is perfectly possible that
$X^{\T}$ is empty, and then it is not clear how to proceed.  For
example, since 
\[
     \EP{( B (\Z/N))} \cong \fmlgpof{C}[N], 
\]
one expects that 
\[
   \Egkv{\T}{\T/\T[N]}  \cong \Egkv{\T[N]}{\point} \cong C[N].
\]
We shall cite three facts to support our conjecture.  
First, if
$F\subseteq X^{\T}$ is a component of the fixed set, and if 
\[
m  : \T\to T
\]
is a reduction of the action of $\T$ on $F$, then 
using the
map of diagram \eqref{eq-extend-equiv-cc-over-cq}
\[
    \ET{\PB (m)} = \EP{\PB (m)} + m : \ET{F}  \rightarrow \VG{\CQ}/W
(m) \rightarrow  \VG{\CQ}/W,
\]
it is easy to check that 
\[
    \anomaly{V}\restr{\ET{X^{\T}}} \cong \anomaly{V'}\restr{\ET{X^{\T}}}
\]
when 
\[
   c_{2} (\borel{V}) \cong c_{2} (\borel{V'}):
\]
this is the content of the proof of Lemma \ref{t-le-gamma-B}.

Second, the isomorphism 
\[
       \ETS{V} \cong \I (V)
\]
of the conjecture corresponds to the fact that the sigma function
gives a trivialization of $\ETS{V}$ in Theorem
\ref{t-th-sigma-orientation}.

Finally, we construct a map $\ED{\PB}$
for a stylized functor $\EDsym$, which is not quite Grojnowski's elliptic
cohomology, but captures its behavior on stalks.  We simply throw away
the points of $C$ for which we have no instructions for  constructing a
map \eqref{eq:15}. The functor $\EDsym$
is inspired by the \emph{rational} $\T$-equivariant elliptic spectra
of Greenlees \cite{Greenlees:Ell} and by Hopkins's study of characters
and elliptic cohomology \cite{Hopkins:Ell}.  

Recall that $\ringedspaces$ denotes the category of ringed spaces.
Let $\OverS$ be the category in which the objects are ringed spaces
$(X,\O{X})$, and in which a map 
\[
   f = (f_{1},f_{2}): (X,\O{X}) \rightarrow (Y,\O{Y})
\]
is a map of spaces
\[
     f_{1}: X\to Y
\]
and a map of sheaves of algebras over $X$
\[
    f_{2}: \O{X}\rightarrow f_{1}^{-1} \O{Y}.
\]
Let $\SubT$ be the category of closed subgroups of $\T$ with
morphisms given by inclusions.   We shall define $\EDsym$ to be a functor 
\[
   \EDsym: \CategoryOf{$\T$-spaces} \rightarrow 
                 \OverS^{\SubT} 
\]
from $\T$-spaces to the category of $\SubT$-diagrams in $\OverS$.

Let $X$ be a $\T$-space.  The ringed space $\ED{X} (\T)$ is just 
\[
     \ED{X} (\T) = \E{X^{\T}}\times C
\]
(which is empty if $X^{\T}$ is).
If $\T[N]\subset \T$ is a finite subgroup and $X^{\T[n]}$ is empty,
then 
\[
     \ED{X} (\T[N])=\emptyset.
\]
Otherwise, 
\[
     \ED{X} (\T[N]) = C[N],
\]
with structure sheaf $j^{-1} \ETS{X}$, where $j$ denotes the
inclusion 
\[
   j: C[N]\to C.
\]
Explicitly, for $U\subseteq C[N]$, $\ED{X} (\T[N]) (U)$ is the product
\[
      \ED{X} (\T[N]) (U) = 
\prod_{a\in U} \ETS{X}_{a} = \prod_{a\in U} \ETS{X^{A}}_{a},
\]
over $a\in U$ of the stalks of Grojnowski's elliptic cohomology.  For
$A=\T[N]\subseteq B$ with $X^{B}$ empty, the map 
\[
    \ED{X} (A) \rightarrow \ED{X} (B)
\]
is trivial; otherwise it is induced by the map of sheaves of algebras
over $C$ 
\[
    \ETS{(X^{A})} \rightarrow \ETS{(X^{B})}
\]
by restriction to $C[N]$.

\begin{Proposition} \label{t-pr-principle-ed}
Let $G$ be a simple and simply connected Lie group, and let $\PB$ be a
principal $G$-bundle over a $\T$-space $X$.  Then there
is a canonical map 
\[
     \ED{\PB} : \ED{X} \rightarrow \VG{C}/W
\]
in $\OverS^{\SubT}$, where $\VG{C}/W$ is considered as a constant
$\SubT$-diagram, such that the diagram 
\[
\begin{CD}
   \TEP{X} @> \TEP{\PB} >> \VG{\fmlgpof{\CQ}}/W \\
   @VVV @VVV \\
   \ED{X} (0) @> \ED{\PB} (0) >> \VG{\CQ}/W
\end{CD}
\]
commutes, and such that for all components $F$ of $X^{\T}$ and all reductions 
\[
m: \T\to T
\]
of the action of $\T$ on $\PB\restr{F}$, 
\[
    \ED{(\PB\restr{F})} (\T) = \E{\PB (m)} + m.
\]
\end{Proposition}

The proof will occupy the rest of this section.  

Let $A$ be a closed
subgroup of $\T$. Let $Y$ be a
connected component of $X^{A}$.  A reduction 
\[
m: A\rightarrow T
\]
of the action of $A$ on $\PB\restr{Y}$ determines a map 
\[
   C (A) \xrightarrow{C\otimes m} (\VG{\CQ})^{W (m)};
\]
as usual if $a\in C (A)$ we write $a^{m}$ for the resulting element of
$(\VG{\CQ})^{W (m)}$.  The
addition  
\[
    \VG{\CQ}\times \VG{\CQ} \xrightarrow{+} \VG{\CQ}
\]
induces a translation 
\[
    (\VG{\CQ})^{W (m)} \times \VG{\CQ}/W (m) \xrightarrow{+} \VG{\CQ}/W (m), 
\]
and so we get an operator
\[
   \trans_{a^{m}} : \VG{\CQ}/W (m) \rightarrow \VG{\CQ}/W (m).
\]

If $A=\T$,
then we define 
\[
   \ED{\PB} (\T)_{Y,m}  = \E{\PB (m)} + m,
\]
as required by the Proposition.  If $A$ is finite,
then we define the map of ringed spaces
$\ED{\PB} (A)_{Y,m,a}$ to be the composition 
\[
\xymatrix{
{\ED{Y} (A)_{a}}
 \ar[rr]^-{\ED{\PB} (A)_{Y,m,a}}
 \ar[d]_{\trans_{-a}}
& &
{\VG{\CQ}/W (m)}
\\
{\ED{Y} (A)_{0}} 
 \ar[r]^{\cong}
& 
{Y_{\HHsym,0}}
 \ar[r]^-{\stalk{(\PB_{\HHsym})}}
&
{\VG{\CQ}/W (m).}
 \ar[u]_{\trans_{a^{m}}}
}
\]
Here we have used the fact that the $\ED{Y} (A)_{0}$ is just the
origin in $C$, with ring 
\[
\stalk{(\O{\ET{Y}})} \cong \HHstalk{Y}; 
\]
Lemma \ref{t-le-hol-char-classes-I} provides the map $\stalk{(\PB_{\HHsym})}$.  We define
\[
     \ED{\PB} (A)_{Y,m} = \coprod_{a\in C[N]} \ED{\PB} (A)_{Y,m,a}: 
   \ED{Y} (A) \rightarrow \VG{\CQ}/W (m).
\]

\begin{Lemma}  \label{t-le-QETYmm}
If $m$ and $m'$ are two reductions of the action of
$A$ on $\PB\restr{Y}$, then the diagram 
\begin{equation} \label{eq:24}
\begin{CD}
\ED{Y} (A) @> \ED{\PB} (A)_{Y,m}  >> \VG{\CQ}/W (m) \\
@V \ED{\PB} (A)_{Y,m'} VV             @VVV \\
\VG{\CQ}/W (m') @>>> \VG{\CQ}/W
\end{CD}
\end{equation}
commutes.
\end{Lemma}

\begin{proof}
This follows from the fact, proved in Lemma
\ref{t-le-g-in-NT-if-Z-conn}, that $m$ and $m'$ differ by an element
of the Weyl group of $G$.
\end{proof}

The lemma permits us to write $\ED{\PB} (A)_{Y}$ for the map
\[
    \ED{Y} (A) \rightarrow \VG{\CQ}/W
\]
described by~\eqref{eq:24}.  We define
\[
\ED{\PB} (A):   \ED{X} (A) = \coprod_{Y} \ED{Y} (A)
    \xra{\coprod \ED{\PB} (A)_{Y}} 
    \VG{\CQ}/W,
\]
where the coproduct is over the components $Y$ of $X^{A}$.  The maps
$\ED{\PB} (A)$ as $A$ ranges over closed subgroups of $\T$ assemble to
give the map $\ED{\PB}$ of Proposition \ref{t-pr-principle-ed}.

\subsection{Relationship to the theory $\ETsym$}

The construction in Proposition \ref{t-pr-principle-ed} is closely
related to the theory $\ETsym$.  We shall briefly explain how
this works, as the explanation sheds light on the
relationship between the ``transfer argument'' of Bott-Taubes and the
geometry of the variety $\VGW{C}.$

Suppose that $a$ is a special point of $C$ of order $N$ and let
$A=\T[N]$.  By definition we have a map 
\[
    \ED{X} (A)_{a} \rightarrow \ETa{X}{a};
\]
indeed it is the inclusion of the stalk at $a$.
Suppose that $b$ is an ordinary point, and let $U = U_{a}\cap 
U_{b}$.  Suppose that $F$ is a component of $Y^{\T}$.  Let
\[
m_{F}: \T\to T
\]
be a reduction of the action of $\T$ on $\PB\restr{F}$; we write
\[
m_{Y} = m_{F}\restr{A}
\]
for the resulting reduction of the action of $A$ on $\PB\restr{Y}$
as in Lemma \ref{t-le-m-Y-restr-to-Fcheck}.\nl

Consider the diagram: 
\begin{equation}\small \label{eq-trans-form-diagram}
\xymatrix{
{\ETa{X}{a}}
&
{\ED{Y} (A)_{a}}
 \ar[rrr]^-{\ED{\PB} (A)_{Y,m_{Y},a}}
 \ar[d]_{\trans_{-a}}
 \ar[l]
&
&
&
{\VG{\CQ}/W (m_{Y})}
 \ar `r[ddr]
      [ddr]
\\
{\ETa{X}{a}\restr{U}} 
 \ar@{^{(}->}[u]
 \ar[dd]_{\glue_{ab}}
&
{\ED{Y} (A)_{0}}
 \ar[rrr]^-{\stalk{(\PB_{\HHsym})}}
&
&
&
{\VG{\CQ}/W (m_{Y})}
 \ar[u]_{\trans_{a^{m_{Y}}}}
\\
&
{\ED{F} (\T)_{0}}
 \ar[rrr]^{\ED{\PB} (\T)_{0}}
 \ar[u]
 \ar[d]
&
&
&
{\VG{\CQ}/W (m_{F})}
 \ar@{->>}[u] 
&
{\VG{\CQ}/W}
\\
{\ETa{X}{b}\restr{U}}
 \ar@{_{(}->}[d]
&
{\ED{F} (\T)}
 \ar[rrr]^-{\ED{\PB} (\T)}
&
&
&
{\VG{\CQ}/W (m_{F})}
 \ar[u]
\\
{\ETa{X}{b}}
&
{\E{F}\times C}
 \ar[rrr]^-{\E{\PB (m_{F})} + m_{F}}
 \ar[u]^{\trans_{-a}}
 \ar[l]
&
&
&
{\VG{\CQ}/W (m_{F})}
 \ar[u]_{\trans_{(-a)^{m_{F}}}}
 \ar `r[uur]
       [uur]
}
\end{equation}

The commutativity of the rectangle on the left is just the definition
of $\glue_{ab}$.  The commutativity of the rectangle on the right is
evident.  The commutativity of the top and bottom rectangles in the
middle is the definition of $\ED{\PB}.$  The commutativity of the
remaining rectangles in the middle follows from the group structure on
$\CQ$ and $\VG{\CQ}$, together with the definitions of the maps involved.

\begin{Remark}\label{rem-1} 
We conclude this paper where the research for it began, with an
explanation of the relationship between ``transfer formula'' of
\cite{BottTaubes:Rig} and the diagram \eqref{eq-trans-form-diagram}.  
Let $F\subseteq Y^{\T} \subseteq Y \subseteq X^{\T[N]}$ be as above.  Let 
\[
m \in \cochars = \hom (\T,T)
\]
be a reduction of the action of $\T$ on $\PB\restr{F}$ (so
$m_{Y} = m\restr{\T[N]}$
is a reduction of the action of $\T[N]$ on $\PB\restr{Y}$).
Let $\theta \in \O{} (\VG{\Gaan})^{W}$ be a theta function for $G$; it
determines a holomorphic characteristic class for principal
$G$-bundles of the form $\borel{\PB}$: that is, 
the characteristic class $\theta (\borel{\PB})$ lies in $\HHol{X}{\affan}$.

The first point is that, for any $a\in \Gaan$,  $\trans_{a^{m}}\theta \in
\O{} (\VG{\Gaan})^{W (m)}$, so it gives a holomorphic
characteristic class for principal $Z (m)$-bundles.
Moreover, the commutativity of the diagram 
\[
\begin{CD}
\HC{F}\times \Gaan @> \TH{\PB (m)}    >> \VG{\Gaan}/W (m) \\
@V \E{F}\times \trans_{a} VV @VV \trans_{a^{m}} V \\
\HC{F}\times \Gaan @> \TH{\PB (m)}    >> \VG{\Gaan}/W (m)
\end{CD}
\]
implies that
\[
       \trans_{a} (\theta (\PB\restr{F})) = 
       (\trans_{a^{m_{F}}}\theta) (\PB (m_{F})) \in \HHol{F}{\affan}.
\]
The second point is that, if $a\in C[N]$,  and $\bara$ is a lift of
$a$ to $\affan$, then 
$\trans_{\bara^{m}}\theta$ is nearly invariant under the action of
$W (m_{Y})$.  Precisely, if $w\in W (m_{Y})$, then 
\[
     \bara^{w m} = \bara^{m} + \lambda
\]
for some $\lambda\in \VG{\Lambda}$: that is, $\bara^{w m}$ and $\bara^{m}$ are
related by the action of the \emph{affine Weyl group} of $G$.  Since
$\theta$ is a theta function for $G$, 
the relationship between 
$
\trans_{\bara^{m}}\theta
$
and 
$
\trans_{\bara^{wm}}\theta
$
is controlled by  $c_{2} (\borel{\PB})$
When this class is zero, or when the second Chern class of another
bundle cancels it, then we may suppose that we have a characteristic
class 
\[
 (\trans_{\bara^{m_{Y}}}\theta) (\PB (m_{Y})) \in \HHol{Y}{\affan}.
\]
We then 
have  
\[
       (\trans_{a^{m_{Y}}}\theta) (\PB (m_{Y}))\restr{F} = 
       (\trans_{a^{m}}\theta) (\PB (m)) = 
       \trans_{a} (\theta (\PB\restr{F})),
\]
which is a typical ``transfer formula''.
\end{Remark}

\subsection{The nonequivariant case}

The conjecture is interesting already in the nonequivariant case.  In
order to compare with \cite{AHS:ESWGTC}, we suppose that $V$ is an
$SU (d)$-bundle over a space $X.$  Let $T\subset SU (d)$ be the usual
maximal torus, with Weyl group $W$.  Let $C=\C/\Lambda$ be a complex
elliptic curve, and let $E$ be the associated elliptic spectrum.

We then have a map 
\[
   X \rightarrow BSU (d)
\]
which in $E$-theory gives (by the splitting principle) a map 
\[
   X_{E} \xrightarrow{f} BSU (d)_{E} \cong (\cochars \otimes \fmlgpof{C})/W.
\]
The line bundle $\I (V) = f^{*}\I (\sigma_{d})$ is certainly
canonically isomorphic to $V_{E}$: this follows simply from the fact
that the sigma function is of the form 
\[
   \sigma (z) = z + o (z^{2}).
\]
Of course the line bundle $\anomaly{V}$ is trivial when restricted to
$\VG{\fmlgpof{C}}/W$.  But if $c_{2} V = 0$, then $\anomaly{V}$ has a
canonical trivialization, since $\anomaly{V}$ descends from the line
bundle $\Loo 
(c_{2})$ over $\VG{\CQ}$ defined by 
\[
   \Loo (c_{2})  =  \frac{\VG{(\Gman)} \times \C}
               {(u,\lambda) \sim 
                (uq^{m}, u^{\Ihat{m}}q^{\qfm{m}} \lambda)}
\]
(see \eqref{eq:21}). It follows that $\sigma (V)$ gives a trivialization of 
\[
     \anomaly{V}\otimes \I (V) \cong V_{E}.
\]
Notice that we only needed $c_{2}V = 0$ to get a trivialization of
$\anomaly{V}$: this is because our elliptic curve is of the form $C=\C/\Lambda$,
and the construction of $\Loo (c_{2})$ uses the covering of $C$.

\end{document}